\documentclass[]{interact}

\usepackage{epstopdf}% To incorporate .eps illustrations using PDFLaTeX, etc.
\usepackage[caption=false]{subfig}% Support for small, `sub' figures and tables

\usepackage[numbers,sort&compress]{natbib}% Citation support using natbib.sty
\bibpunct[, ]{[}{]}{,}{n}{,}{,}% Citation support using natbib.sty
\renewcommand\bibfont{\fontsize{10}{12}\selectfont}% Bibliography support using natbib.sty
\makeatletter% @ becomes a letter
\def\NAT@def@citea{\def\@citea{\NAT@separator}}% Suppress spaces between citations using natbib.sty
\makeatother% @ becomes a symbol again
\usepackage{anyfontsize}
\theoremstyle{plain}% Theorem-like structures provided by amsthm.sty
\newtheorem{theorem}{Theorem}[section]
\newtheorem{lemma}[theorem]{Lemma}
\newtheorem{corollary}[theorem]{Corollary}
\newtheorem{proposition}[theorem]{Proposition}

\theoremstyle{definition}
\newtheorem{definition}[theorem]{Definition}
\newtheorem{remark}[theorem]{Remark}
\newtheorem{example}[theorem]{Example}

\usepackage{xcolor}

\theoremstyle{remark}

\begin{document}
	
	%\articletype{ARTICLE TEMPLATE}% Specify the article type or omit as appropriate

	\title{Hyers-Ulam Stability of Unbounded Closable Operators in Hilbert Spaces}
	\author{
		\name{Arup Majumdar\textsuperscript{a},  P. Sam Johnson\textsuperscript{b}, Ram N. Mohapatra\textsuperscript{c}}
		\affil{\textsuperscript{a,}\textsuperscript{b}Department of Mathematical and Computational Sciences,
			National Institute of Technology Karnataka (NITK), Surathkal, Mangaluru 575025, India.} \textsuperscript{c} Department of Mathematics, University of Central Florida, Orlando, FL. 32816, USA}

	\maketitle
	
	\begin{abstract}
		In this paper, we discuss the Hyers-Ulam stability of closable  (unbounded) operators with several interesting examples.  We also 	
		present results pertaining to the Hyers-Ulam stability of the sum and product of closable operators to have the Hyers-Ulam stability and the necessary and sufficient conditions of the Schur complement and the quadratic complement of $2 \times 2$ block matrix $\mathcal A$ in order to have the Hyers-Ulam stability.
	\end{abstract}
	
	\begin{keywords}
		Hyers-Ulam stable operator, closed range operator, closable operator.
	\end{keywords}
      \begin{amscode}47A05, 47A10, 47A55.\end{amscode}

	\section{Introduction}
	The theory of Hyers-Ulam stability plays an important role in functional equations, optimization, differential equations and many branches of mathematics and statistics. Ulam, in a talk delivered at the University of Wisconsin in 1940, discussed some important unsolved problems and  raised the stability problem for functional equations, which reads as: ``For what metric groups $G$, is it true that an $\varepsilon$-automorphism of $G$ is necessarily near to a strict automorphism?" The question was answered by Hyers in 1941 for Banach spaces: Let $X$ and $Y$ be two real Banach spaces and $f: X \to Y$ be a mapping such that for each fixed $x \in X$, $f(tx)$ is continuous in $t\in \mathbb{R}$, the set of all real numbers. If there exists  $\varepsilon \geq 0$ such that $$\| f(x+y) - f(x) -f(y)\| \leq \varepsilon \quad \text{ for all }x, y \in X,$$ then there exists a unique linear mapping $L : X \to Y$ such that $\|f(x) - L(x)\| \leq \varepsilon$ for every $x \in X$. This result is called the Hyers-Ulam stability of the additive Cauchy equation $$g(x+y) = g(x) + g(y).$$
	In 1978, Rassias considered the unbounded right-hand side in the involved inequalities, depending on certain functions, now known as the modified Hyers-Ulam  stability for the additive functional equation 
	\cite{MR1778016}. After that, many researchers have extended Ulam's stability problems to other functional equations and generalized Hyers' result in various directions. Over the last three decades, this topic has been very well known as Hyers-Ulam stability, or sometimes it is referred to as Hyers-Ulam-Rassias stability.  

	Obloza \cite{MR1321558} was the first author who proved results concerning the Hyers-Ulam stability of differential equations. Alsina and Ger \cite{MR1671909} investigated the Hyers-Ulam stability of first-order linear differential equations. Miura et al. generalized the results for $n^\text{th}$ order linear differential operator $p(D)$ and proved that the differential operator equation $$p(D)f=0$$ is Hyers-Ulam stable iff the algebraic equation $p(z)= 0$ has no pure imaginary solution, where $p$ is a complex-valued polynomial of degree $n$, and $D$ is a differential operator \cite{MR2000046}. Moreover, they introduced the notion of the Hyers-Ulam stability of a mapping (not necessarily linear) between two complex linear spaces $X$ and $Y$ with gauge functions $\rho_X$ and $\rho_Y$, respectively. A mapping $f$ has the Hyers-Ulam stability (HUS) if there exists a constant $M \geq 0$ with the following property:
	
	For every $\varepsilon \geq 0$, $y\in f(X)$ and $x\in X$ satisfying $\rho_Y (f(x) -y) \leq \varepsilon$ we can find an $x_0 \in X$ such that $f(x_0) =y$ and $$\rho_X (x-x_0) \leq M \varepsilon,$$ where $M$ is called as a HUS constant, and let us denote the infimum of all the HUS constants for $f$ by $M_f$. In other words, if $f$ has the HUS, then for each $y\in f(X)$ and ``$\varepsilon$-approximate solution" $x$ of the equation $f(u)= y$ there corresponds an exact solution $x_{0}$ of the equation that is contained in a $M\varepsilon$- neighbourhood of $x$.

	From now on, by an operator, we shall mean a non-zero linear operator. We will analyse the Hyers-Ulam stability of unbounded operators between Hilbert spaces. The specification of a domain is an essential part of the definition of an unbounded operator, usually defined on subspaces.  Consequently, for an operator $T$, the specification of the subspace $D$ on which $T$ is defined, called the domain of $T$, denoted by $D(T)$, is to be given. The null space and range space of $T$ are denoted by $N(T)$ and $R(T)$, respectively. $W_{1}^{\perp}$ denotes the orthogonal complement of a set $W_{1}$ whereas $W_{1} \oplus^{\perp} W_{2}$ denotes the orthogonal direct sum of the subspaces $W_{1}$ and $W_{2}$ of $X$. Moreover,  $T_{W}$ denotes the restriction of $T$ to a subspace $W$ of  $D(T)$.  We call $D(T)\cap N(T)^\perp $, the carrier of $T$ and it is denoted by $ C(T)$.  
	For the sake of completeness of exposition, we first begin with the definition of a closed operator.

	\begin{definition}
		Let $ T$ be an operator from a Hilbert space $H $ with domain $D(T) $ to a Hilbert space $K $. If the graph of $T$ is defined by 
		$$ G(T)=\left\{(x,Tx): x\in D(T)\right\} $$ is closed in $H\times K $, then $T$ is called a closed  operator. Equivalently, $T $ is a closed operator if $ \{x_n \}$ in $D(T) $ such that $ x_n\rightarrow x $ and $ Tx_n\rightarrow y$ for some $ x\in  H,y\in   K $,  then $ x\in  D(T) $ and $ Tx=y $. That is, $G(T)$ is  a closed subspace of $H\times K$ with respect to the graph norm $\|(x, y)\|_T=(\|x\|^2+\|y\|^2)^{1/2}$. It is easy to show that the graph norm  $\|(x, y)\|_T$ is equivalent to the norm $\|x\|+\|y\|$. 	 We note that,  for any densely defined closed operator $T$,  the closure of $C(T)$, that is,  $ \overline{C(T)}$ is $ N(T)^\perp$. 
		We say that $S$ is an extension of $ T $ (denoted by $ T\subset S$) if $ D(T)\subset D(S) $ and $ Sx=Tx $ for all $ x\in D(T)$.
		
		An operator $ T $ is said to be closable 	
		if $T$ has a closed extension. It follows that $T$ is closable iff the closure $\overline {G(T)}$ of $G(T)$ is a graph of an operator. 	It is also possible for a closable operator to have many closed extensions.  Its minimal closed extension is denoted by $ \overline{T}$. That is, every closed extension of $ T $ is also an extension of $ \overline{T}$.
		
	\end{definition}

	Let $T$ be a closed operator from $D(T) \subset H$ to $K$. We define the induced one-one operator $\tilde{T}$ from $C(T)$ into $K$ by $
	\tilde{T}x = Tx,  \text{ for  }  x\in C(T)$, which is also closed.

	\begin{definition}
		Let $T$ be a closed operator from $D(T) \subset H$ to $K$. We say that $T$ is lower semibounded if there exists  $\gamma> 0$ such that
		\begin{equation*}
			\|Tx\| \geq \gamma \|x\|,  \quad \text{ for all } x\in C(T).
		\end{equation*}
		The minimum modulus of $T$ is defined by 
		\begin{align*}
			\gamma(T) &= \sup\Big\{\gamma : \|Tx\| \geq \gamma \|x\|, \text{ for all } x\in C(T) \Big\}.
		\end{align*}
		
	\end{definition}
	\begin{theorem}\cite{MR0200692}\label{thm 1}
		Let $T$ be a closed operator from $D(T) \subset H$ to $K$. The range of $T$ is closed iff $\gamma(T) >0$.
	\end{theorem}

	\begin{definition}\cite{MR2953553}
		A linear subspace $D$ of $D(T)$ is called a core for $T$ if $D$ is dense in $(D(T), \|.\|_{T})$, that is, for each $x\in D(T)$, there exists a sequence $\{x_{n}\}$ in $D$ such that  $\{x_{n}\}$ and $\{ Tx_{n}\}$ converge to $x$ and $Tx$ respectively. Moreover, if $T$ is closed, a linear subspace $D$ of $D(T)$ is a core for $T$ iff $T=\overline{T_{D}}$.
	\end{definition} 
	\begin{remark}
		Let $T$ be a closed operator from  $D(T)\subset H$ to  $K$. If $N(T) \neq \{0\}$, then $C(T)$ is not a core for $T$. In this case, 	 $C(T)$ can not be dense in $D(T)$. Indeed, suppose that  $C(T)$ is dense in $D(T)$, then 	$ N(T)\subset D(T)\subset \overline{C(T)}$. Since $N(T) \neq \{0\}$,  there exists $x\in N(T)$ with $x \neq 0$. For each $n\in \mathbb N$, there is an element $x_{n}\in C(T)$ such that $\|x- x_{n}\|_T < \frac{1}{n}$, hence  $\|x\|^{2} + \|x_{n}\|^{2} < \frac{1}{n^{2}}$, thus $\|x\| < \frac{1}{n}$. Since $n$ is arbitrary, we get that $x=0$, which is a contradiction. 
	\end{remark}

	\begin{definition}
		Let $T$ be a closed operator from  $D(T) \subset H$ to  $K$. The generalized inverse of $T$ is the map $T^{\dagger}: R(T) \oplus^{\perp} R(T)^{\perp} \to H$ defined by
		\begin{equation}\label{equ 1}
			T^{\dagger} y = 
			\begin{cases}
				(T_{C(T)})^{-1}y  & \text{if} ~  y\in R(T)\\
				0    & \text{if}  ~ y\in R(T)^{\perp}.
			\end{cases}
		\end{equation}
		
	\end{definition}
	\begin{definition}
		A linear operator $L$ from $D(L) \subset H$ to  $K$ has a decomposable domain if $D(L) = N(L) \oplus^{\perp}  C(L)$, where $H = \overline{N(L)} \oplus^{\perp} {N(L)}^{\perp}$. The generalized inverse of $L$ is the map $L^{\dagger}: R(L) \oplus^{\perp} R(L)^{\perp} \to H$ defined as follows:
		\begin{equation}\label{equ 2}
			L^{\dagger} y = 
			\begin{cases}
				(L_{C(L)})^{-1}y  & \text{if} ~  y\in R(L)\\
				0    & \text{if}  ~ y\in R(L)^{\perp}.
			\end{cases}
		\end{equation}
	\end{definition}

	\noindent We now see the definition of Hyers-Ulam stability for linear operators between normed spaces.   
	\begin{definition}
		Let $T$ be a linear operator from $D(T) \subset X$ to $Y$, where $X$ and $Y$ are both normed linear spaces. The operator $T$ is said to have the Hyers-Ulam stability if there exists a constant $M \geq 0$ with the following equivalent properties : 
		\begin{enumerate}
			\item For any $y\in R(T), ~ \varepsilon \geq 0$ and $x\in D(T)$ with $\| Tx -y\| \leq \varepsilon$, there exists  $x_0 \in D(T)$ such that $Tx_{0} = y$ and $\|x -x_{0}\| \leq M\varepsilon$.
			
			\item For any $\varepsilon \geq 0$, $x\in D(T)$ with $\| Tx\| \leq \varepsilon$, there exists $x_0 \in D(T)$ such that $Tx_{0} =0$ and $\|x -x_{0}\| \leq M \varepsilon$.
			\item For any $x\in D(T)$, there exists $x_0 \in N(T)$ such that $\|x-x_{0}\| \leq M \|Tx\|$.

		\end{enumerate}
		We call $M$ a HUS constant for $T$.  The infimum of all HUS constants for $T$ is denoted by $M_T$, and it need not be a HUS constant.  It is proved in \cite{MR2111613} that $M_T$ is a HUS constant for $T$ when $N(T)$ is proximinal. We also say that $M_T$ is a HUS constant when $T$ is considered as a closed operator between Hilbert spaces.
		
We characterise closable Hyers-Ulam stable operators, their closures and generalized inverses in Section 2.  We investigate the relationship between a densely defined closable Hyers-Ulam stable operator $T$ with the spectrum of ${T^{*}\overline{T}}_{N(T^{*}\overline T)^{\perp}}$ and we analyze Hyers-Ulam stability of  $T^{n}$, $n\in \mathbb N$,  for the paranormal operator $T$. We also discuss the sum and product of Hyers-Ulam stable operators to be again Hyers-Ulam stable and the necessary and sufficient conditions of the Schur complement and the quadratic complement of $2 \times 2$ block matrix $\mathcal A$ in order to have the Hyers-Ulam stability in Section 3. In the final section, we illustrate an example of a densely defined closable Hyers-Ulam stable operator which is not closed. 
		
	\end{definition}

	\section{Characterizations of Hyers-Ulam Stable Operators}
	We start this section with two examples. We consider two operators, namely, the Bernstein operator and the classical $n$-th Sz{\'a}sz-Mirakjan operator \cite{MR3753562}. We denote the space of all continuous real-valued functions defined on the interval $[0,1]$ by $\mathcal C[0,1]$ whereas the space of all bounded real-valued continuous functions on $[0, +\infty)$ by $\mathcal C_{b}[0, +\infty)$. We examine the HUS of Bernstein and  $n$-th Sz{\'a}sz-Mirakjan operators using the definition of HUS for linear operators.

	\begin{example}
		The Bernstein operator $B_{n}: (\mathcal{C}[0,1]  , \|.\|_{\infty}) \to  (\mathcal{C}[0,1]  , \|.\|_{\infty})$ is defined by $$(B_{n}f)(x) = \sum_{k=0}^{n} f\Big(\frac{k}{n}\Big)P_{n,k}, \quad \text { for all } x\in [0,1],$$ where $P_{n,k}=
		\binom{n}{k} x^{k} (1-x)^{n-k}$ and $\|.\|_{\infty}$ is defined as sup-norm. 
		 Then $\|B_{n} f\|_{\infty} \leq \|f\|_{\infty},$ for all  $f\in \mathcal{C} [0,1]$, hence $B_{n}$ is a bounded operator. Since  $\{P_{n.k}\}_{k=0}^{n}$ is linearly independent, the dimension of  $R(B_{n})$ is $n+1$. The induced one-one operator $\tilde{B_n}$ is bounded and its inverse $(\tilde{B_{n}})^{-1}$ exists because $R(B_{n})$ is finite dimensional.
		 
		  For given $\varepsilon > 0$,  $B_{n} g = \sum_{k=0}^{n} g(\frac{k}{n}) P_{n,k} \in R(B_{n})\setminus \{0\}$ with $\|B_{n}g\|_{\infty} < \varepsilon$, there exists an element $h\in \mathcal{C}[0,1]$ such that $h(\frac{k}{n}) = g(\frac{k}{n})  ,~ 0\leq k \leq n$, and $h(x)$ is the straight line on $[\frac{k}{n}, \frac{k+1}{n}]$ by joining two points $(\frac{k}{n}, g(\frac{k}{n}))$ and $(\frac{k+1}{n}, g(\frac{k+1}{n}))$, for all $0\leq k \leq (n-1)$. Therefore $h\notin N(B_{n})$, hence 	  
		   $(\tilde{B_{n}}) (h + N(B_{n})) = B_{n}(h) = B_{n}(g)$ and $(\tilde{B_{n}})^{-1} (B_{n}(g)) = h+ N(B_{n})$.
		   We have $\|h + N(B_{n})\| \leq \|(\tilde B_{n})^{-1}\| \|B_{n} (g)\|_{\infty} < \|(\tilde B_{n})^{-1}\| \varepsilon$, so there exists    
		    $t\in N(B_{n})$ such that $\|h + t\|_{\infty} \leq \|(\tilde B_{n})^{-1}\| \varepsilon$. Hence  $$\|g-g_0\|_{\infty} \leq \|(\tilde B_{n})^{-1}\| \varepsilon, \quad \text{ where }g_0=g-h-t \in N(B_{n}).$$ Therefore, $B_{n}$ is Hyers-Ulam stable.
	\end{example}
	\begin{example} Let us consider
	$n$-th Sz{\'a}sz-Mirakjan operator $M_{n}: (\mathcal C_{b}[0, +\infty), \|.\|_{\infty}) \to (\mathcal C_{b}[0, +\infty), \|.\|_{\infty})$ defined by $$M_{n}f(x) = e^{-nx} \sum_{k=0}^{\infty} f\Big(\frac{k}{n}\Big) {\frac{n^{k} x^{k}}{k!}}.$$ We claim that $M_{n}$ is not  Hyers-Ulam stable.
		
		Suppose that $M_{n}$ is Hyers-Ulam stable. Then there exists $N>0$ with the following property:  for each $f\in \mathcal C_{b}[0, +\infty)$ with $\|M_{n}f\|_{\infty} \leq 1$, there exists $g \in N(M_{n})$ with $\|f-g\|_{\infty} \leq N$. 
		
		According to the Stirling formula, we have $\frac{i^{i}}{i! e^{i}} \sim \frac{1}{\sqrt{2\pi i}}$, hence $\displaystyle \lim_{i\to \infty} \frac{i^{i}}{i ! e^{i}}  =0$, so there exists  $j\in \mathbb{N}$ such that $(N+1) \frac{j^{j}}{j ! e^{j}} \leq 1$. 	
We define $f$ by 	
		$$ f(x)= 	\begin{cases}
		0,     & \text{if} ~ x\in [0, \frac{j-1}{n}] \cup [\frac{j+1}{n}, \infty)\\
		n(N+1)x -(N+1)(j-1) & \text{if} ~ x\in [\frac{j-1}{n}, \frac{j}{n}]\\
		-(N+1)nx + (N+1)(j+1)& \text{if} ~ x\in [\frac{j}{n}, \frac{j+1}{n}].
		
	\end{cases} 
	$$
	Hence 	$f\in \mathcal{C}_{b}[0, +\infty)$, 
	$\|f\|_{\infty} = N+1$ and $\|M_{n} f\|_{\infty} =| f(\frac{j}{n} ) \frac{j^{j}}{j! e^{j}}| \leq 1$. Then there exists $g\in N(M_{n})$ such that $\|f-g\|_{\infty} \leq N$. Hence $|f(\frac{j}{n}) -g(\frac{j}{n})| = N+1$, which is a contradiction. Thus, $M_{n}$ is not Hyers-Ulam stable.
	\end{example}

	The Hyers-Ulam stability of a  closed linear operator can be characterized by its closed range, as shown below in the following result.  We shall use this characterization and discuss the HUS of unbounded operators on sequence spaces.
	\begin{theorem}\cite{MR2204863}\label{thm 3}
		Let $T$ be a closed operator from $D(T) \subset H$ into $K$ and $\tilde{T}$ be the induced one-one operator. The following assertions are equivalent:
		\begin{enumerate}
			\item $T$ has the Hyers-Ulam stability.
			\item $T$ has closed range.
			\item ${\tilde{T}}^{-1}$ is bounded.
			\item $T$ is lower semibounded.
		\end{enumerate}
		Moreover, if one of the conditions above is true, then $M_T = \|{\tilde{T}}^{-1}\| = \gamma(T)^{-1}$.
	\end{theorem}
	\begin{example}
		Define $T$ on $\ell^{2}$ by $$T(x_{1}, x_{2}, x_{3},\ldots,x_{n},\ldots) =(x_{1}, 2x_{2}, 3x_{3},\ldots,nx_{n},\ldots )$$ with domain $D(T)=\{(x_{1}, x_{2}, x_{3},\ldots,x_{n},\ldots)\in \ell^{2} : \sum_{n=1}^\infty |n x_{n}|^2<\infty\}$. Since  
		$D(T)$ contains the space $c_{00}$ of all finitely non-zero sequences, $ D(T)$ is a proper dense subspace of $\ell_2$. 
One can show that $T$ is a self-adjoint operator and  ${R(T)} = N(T^{*})^{\perp}= \ell^{2}$. Therefore, $T$ is Hyers-Ulam stable. 

Since $T$ is a densely defined closed operator and injective, $T^{-1}$ exists and is closed.  By the closed graph theorem,  $T^{-1}$ is bounded but $R(T^{-1})$ is not closed. Thus, $T^{-1}$ is not Hyers-Ulam stable. 
	\end{example}
	\begin{example}
		Define $T$ on $\ell^{2}$ by $$T(x_{1}, x_{2},\ldots,x_{n},\ldots)= (x_{2}, 0, 2x_{4},0, 3x_{6},\ldots,0,nx_{2n},0,\ldots)$$ with domain $D(T) = \{(x_{1}, x_{2},\ldots,x_{n},\ldots)\in \ell^{2} : \sum_{n=1}^\infty |n x_{2n}|^2<\infty	 \}$. By sequential argument, one can show that $G(T) = \overline {G(T)}$; hence $T$ is closed. Now we consider operator $\tilde{T}: D(T)/ N(T) \to R(T)$ by $\tilde{T} (x+N(T))= T(x)$, for all $x\in D(T) $. Since $\tilde{T}$ is closed,  $ (\tilde{T})^{-1}$ is also closed. 
		For each $y\in R(T)$,  $\|(\tilde T^{-1}) y\|^{2} \leq \sum_{n=1}^{\infty} {\frac{|y_{2n-1}|^{2}}{n^{2}}} \leq \|y\|^{2}$, so $(\tilde T)^{-1}$ is bounded. By the closed graph theorem, $D((\tilde T)^{-1}) = R(T)$ is closed. Therefore, $T$ is Hyers-Ulam stable.
	\end{example}
	\begin{example}
		Define $T$ on  $\ell^{2}$ by $$T(x_{1}, x_{2},\ldots,x_{n},\ldots)= \Big(x_{1}, 2x_{2}, \frac{x_{3}}{3}, 4x_{4}, \frac{x_{5}}{5},\ldots\Big)$$ with domain $D(T)= \{(x_{1}, x_{2}, x_{3},\ldots,x_{n},\ldots)  : (x_{1}, 2x_{2}, \frac{x_{3}}{3}, 4x_{4}, \frac{x_{5}}{5},\ldots)\in \ell^{2} \}$. One can show that $T$ is closed and $R(T)$ is a proper dense subspace of $\ell_2$.  Therefore, $R(T)$ is not closed, which confirms that $T$ is not Hyers-Ulam stable.
	\end{example}
	\begin{proposition}\label{pro 2}
		Let $T$ be a closable operator from  $D(T) \subset H$ to  $K$. If $T$ is Hyers-Ulam stable, then $\overline{T}$ is  Hyers-Ulam stable with $M_{\overline T} \leq M_{T}$.
	\end{proposition}
	\begin{proof}
		Since $T$ is Hyers-Ulam stable, there exists  $M>0$ with the following property: for any $\varepsilon \geq 0,  u\in D(T)$ with $\|Tu\| \leq \varepsilon$, there exists  $u_{0} \in N(T)$ such that $\|u-u_{0}\| \leq M \varepsilon$.
		
		 Let $x\in D(\overline T)$ with $\|\overline Tx\| \leq {\varepsilon}$. There is a $n_{0}\in \mathbb N$ such that $x_{n_{0}} \in D(T)$ and $\|x_{n_{0}}- x\| \leq \frac{\varepsilon}{m}$ and  $\| Tx_{n_{0}} - \overline{T}x\| \leq \frac{\varepsilon}{m}$, for some large natural  number $m$. Then $\|Tx_{n_{0}}\| \leq \|Tx_{n_{0}} - \overline {T}x\| + \|\overline{T}x\| \leq \frac{m+1}{m} \varepsilon$.
		From the definition of Hyers-Ulam stability of $T$, we get an element in $N(T)$, say $x_{0}$, such that $$\|x_{n_{0}} - x_{0}\| \leq \frac{m+1}{m} M\varepsilon.$$ Hence  $\|x_{0} - x\| \leq \|x_{0} -x_{n_{0}}\| + \|x_{n_{0}}- x\| \leq (\frac{M+1}{m}  + M)\varepsilon$, where $x_{0}\in N(T) \subset N(\overline T)$. Thus, $\overline{T}$ is Hyers-Ulam stable with $M_{\overline T} \leq (\frac{M+1}{m}  + M) $. Here, $m$ is arbitrarily large; we get that  $M_{\overline T} \leq M_{T}$.
	\end{proof}

	\begin{theorem}\label{thm 2}
		Let $T$ be a closed operator from 
		$D(T) \subset H$ to 
 $K$. If $T$ is  Hyers-Ulam stable and $D$ is a core for $T$ with $N(T) = \overline{N(T_{D})}$, then $T_{D}$ is also Hyers-Ulam stable and $M_{T_{D}} \leq M_{T}$.
	\end{theorem}
	\begin{proof}
		From the definition of the Hyers-Ulam stability of $T$, we get  $M_{T}\geq 0$ with the following property: for any $u\in D(T)$, there exists  $u_{0} \in N(T)$ such that $\|u-u_{0}\| \leq M_{T} \|Tu\| < (M_{T}+ \varepsilon) \|Tu\|$, for all $ \varepsilon > 0$.
		 Since $D$ is a core for  $T$, $\overline{T_{D}} = T$.  
		Now, we will show that $T_{D}$ is Hyers-Ulam stable.
		For $u_{0} \in N(T)= \overline {N(T_{D})}$,  there exists a sequence $\{u_{n}\}$ from $N(T_{D})$ such that $u_{n} \to u_{0}$ as $n \to \infty$. We get a natural number, say m, such that $\|u_{m} -u_{0}\| \leq (M_{T} + \varepsilon )\|Tu\| - \|u- u_{0}\|$, where $u_{m} \in N(T_{D})$.
		\begin{equation}\label{ineq-3}
			\|u_{m} -u\| \leq \|u_{m} -u_{0}\| + \|u- u_{0}\| \leq (M_{T} + \varepsilon ) \|Tu\|.
		\end{equation}
		The inequality (\ref{ineq-3}) guarantees that $T_{D}$ is Hyers-Ulam stable and $M_{T_{D}} \leq M_{T}$.
	\end{proof}
 \begin{corollary}\label{cor 1}
 Let $T$ be a closable operator from 
 $D(T) \subset H$ to 
 $K$ with $N(\overline T) = \overline {N(T)}$. Then the following statements hold good :
 \begin{enumerate}
 	\item  
 $T$ is  Hyers-Ulam stable iff $\overline T$ is Hyers-Ulam stable. In this case, $M_{T} = M_{\overline T}$.
 \item If $R(T)$ has a finite deficiency in $K$, then $T$ is Hyers-Ulam stable. 
 \item If $T$ is densely defined, then 
 $T$ is Hyers-Ulam stable iff $T^{*}$ is Hyers-Ulam stable.
						\item If $T$ is densely defined, $T$ is  Hyers-Ulam stable iff $R(T^{*}\overline T)$ is closed.

		\end{enumerate} 
	\end{corollary}
	\begin{proof}
\begin{enumerate}
	\item The proof is obvious from Proposition \ref{pro 2}, and Theorem \ref{thm 2} .
	\item  Since dim$R(T)^{\perp}$= dim$R(\overline T)^{\perp}$ is finite, $R(\overline T)$ is closed.	Hence, $\overline T$ is Hyers-Ulam stable with $N(\overline T)= \overline{N(T)}$, which implies $T$ is Hyers-Ulam stable.
		\item The proof follows from the implication : $R(\overline T)$ is closed iff $R(T^{*})$ is closed.
		\item The proof follows from the implication : $R(\overline T)$ is closed iff $R(T^{*}\overline T)$ is closed.
	
\end{enumerate}	
	\end{proof}

	\begin{lemma}\label{lemma 2}
		Let $T$ be a closable operator from 
		$D(T) \subset H$ to 
		$K$ with an assumption that  $D(T)$ is the decomposable domain.  If $T^{\dagger}$ is Hyers-Ulam stable, then ${\overline T}^{\dagger}$ is  Hyers-Ulam stable.
	\end{lemma}
	\begin{proof}
		First we will show that $R(T)^{\perp} = R(\overline T)^{\perp}$. It is obvious that $ R(\overline T)^{\perp} \subset R(T)^{\perp}$. Let $y\in R(T)^{\perp}$. Then $\langle Tx, y \rangle = 0$, for all $x\in D(T)$.  For each $x\in D(\overline T)$, there exists a sequence $\{x_{n}\}$ in $ D(T)$ such that $\{x_{n}\} $ and $\{Tx_{n}\} $ converge to $x$ and ${\overline T}x$ respectively.
		So 	 $\displaystyle\langle {\overline T}x , y \rangle = \lim_{n\to \infty} \langle Tx_{n}, y\rangle = 0$, therefore $\langle {\overline T}x , y \rangle = 0$, for all $x\in D(\overline T)$ which  implies that $y\in R(\overline T)^{\perp}$. Hence  $R(T)^{\perp} = R(\overline T)^{\perp}$. So, $D(T^{\dagger}) = R(T) \oplus^{\perp} R(T)^{\perp}$ and $D(\overline T^{\dagger}) = R(\overline T) \oplus^{\perp} R(T)^{\perp}$. Thus $N(T^{\dagger}) = N(\overline T)^{\dagger} = R(T)^{\perp}$.
		
		Since $T^{\dagger}$ is Hyer-Ulam stable, there exists  $M > 0$ with the following property: for each $\varepsilon >0$ and $u\in C(T)$ with $\| T^{\dagger} Tu\| =\|u\| \leq \varepsilon$, there exists   $v\in N(T^{\dagger})$ and $\|Tu -v\| \leq M \varepsilon$. Now taking, $(\overline T^{\dagger}) y\in R(\overline T^{\dagger}) = C(\overline T)$, where $y \in R(\overline T)$, and $\|(\overline T^{\dagger}) y\| \leq \frac{\varepsilon}{2}$ with $(\overline T^{\dagger})y =z\in C(\overline T)$ implies $\|z\| \leq \frac{\varepsilon}{2} $ and ${\overline T}z =y$. There exists a sequence $\{z_{n}\}$ of element of $C(T)$ such that $z_{n}\to z$ and $Tz_{n} \to {\overline T}z$ as $n \to \infty$. $T^{\dagger} Tz_{n} = z_{n} $, for all $n\in \mathbb {N}$. We get a natural number $k\in \mathbb N$ such that $\|z_{k}\| < \varepsilon$ and $\| Tz_{k} - {\overline T} z\| < \varepsilon$. Since $T^{\dagger}$ is Hyers-Ulam stable,  there exists $v_{k} \in R(T)^{\perp}$ such that $\|Tz_{k} - v_{k}\| \leq M \varepsilon$. 
		\begin{equation}\label{inequ 4}
			\|y - v_{k} \| =\| {\overline T}z -v_{k}\| \leq \|{\overline T} z - Tz_{k}\| + \|Tz_{k} -v_{k}\| < (M+ 1) \varepsilon.
		\end{equation}
  
  Now we consider an element $w\in D(\overline T^{\dagger})$ with $\| (\overline{T}^{\dagger}) w\| \leq \frac{\varepsilon}{2}$. The inequality (\ref{inequ 4}) confirms that the existence of an element $w_{0} \in R(T)^{\perp} $ such that $\|w-w_{0}\| < (M+1) \varepsilon$. Therefore, ${\overline T}^{\dagger}$ is Hyers-Ulam stable.
	\end{proof}
	\begin{theorem}\label{thm 6}
		Let $T$ be a closable operator from 
		$D(T) \subset H$ to 
		$K$ with an assumption that $D(T)$ is the decomposable domain. If $T^{\dagger}$ is  Hyers-Ulam stable, then $T$ is bounded. 
	\end{theorem}
	\begin{proof}
		As $\overline T$ is closed,  ${\overline T}^{\dagger}$ is closed, where $D({\overline T}^{\dagger}) = R(\overline T) \oplus^{\perp} R(T)^{\perp}$ and $R({\overline T}^{\dagger}) = C(\overline T)$.
		First, we claim that $T$ is bounded when $T^{\dagger}$ is Hyers-Ulam stable.	
Since $T^{\dagger}$ is  Hyers-Ulam stable, by the Lemma \ref{lemma 2}, we get that 
	$R(\overline T^{\dagger}) = C(\overline T)$ is closed in $H$. Also, $N(\overline T)$ is closed.  	
	Let $x\in \overline{ N(\overline T) + C(\overline T)} = \overline{D(\overline T)}$.  Then there exists a sequence $\{x_{n}\}$ in $ D(\overline T)$ such that $x_{n} \to x$ as $n\to \infty$.  Now, we can write $x_{n} = y_{n} + z_{n}$,  for all $n\in \mathbb N$, where $y_{n} \in  N(\overline T)$ and $z_{n} \in C(\overline T)$. So, $\{x_{n}\}$ is Cauchy. Thus $\{y_{n}\}$ and $\{z_{n}\}$ both are Cauchy sequences in $N(\overline T)$ and $C(\overline T)$ respectively. 
		
		As $N(\overline T)$ and $C(\overline T)$ are closed,  $\{y_{n}\}$ and $\{z_{n}\}$  converge to $y$ and $z$ respectively,  for some $y\in N(\overline T)$ and $z\in C(\overline T)$. Thus $x_{n} \to (y+z) \in D(\overline T)$ as $n\to \infty$. So $x= y+z \in D(\overline T)$. Thus, $D(\overline T)$ is closed. Therefore, by the closed graph theorem,  $\overline T$ is bounded; hence $T$ is bounded.

	\end{proof}
	\begin{corollary}
		Let $T$ be a closed operator  from 
		$D(T) \subset H$ to 
		$K$. Then $T$ is bounded iff $T^{\dagger}$ is Hyers-Ulam stable. Moreover, $T$ is bounded iff $T_{C(T)}^{-1}$ is Hyers-Ulam stable with $M_{T_{C(T)}^{-1}} = \|T_{C(T)})\|$.
	\end{corollary}
\begin{proof}
	Suppose that $T$ is bounded.  Then $\overline T$ is bounded, hence $\overline T$ is closed. So  $D(\overline T)$ is closed.  Hence $D(\overline T) \cap N(\overline T)^{\perp}=C(\overline T) = R({\overline T}^{\dagger}) $ is a closed subspace of $H$. Thus  ${\overline T}^{\dagger} = T^{\dagger}$ is Hyers-Ulam stable.  Converse follows from the Theorem \ref{thm 6}. Furthermore, the additional part can be shown from the Theorem \ref{thm 3}.
\end{proof}
	\begin{corollary}
		Let $T$ be a densely defined closable operator  from 
		$D(T) \subset H$ to 
		$K$ with $N(\overline T)= \overline{N(T)}$. Then $|\overline T|$  is Hyers-Ulam stable iff $T$ is a Hyers-Ulam stable opeartor, where $|\overline T| = (T^{*} \overline T)^{\frac{1}{2}}$.
	\end{corollary}
	\begin{proof}
		By the polar decomposition of $\overline T$, we write $\overline T= U_{\overline T} |\overline T|$, where $U_{\overline T}$ is a partial isometry with initial space $N(\overline T)^{\perp} = \overline {R(T^{*}) } = \overline{ R(|\overline T|)}$ and final space $N(T^{*})^{\perp} = \overline{R(\overline T)}$ which gives the equality $R(|\overline T|) = R(T^{*})$. 
As $|\overline T|$ is self-adjoint, 		 $|\overline T|$ is closed. The proof is completed from the implication : 
		$|\overline T|$ is Hyers-Ulam stable iff $R(|\overline T|) = R(T^{*})$ is closed iff $R(\overline T)$ is closed iff $\overline T$ is Hyers-Ulam stable iff $T$ is Hyers-Ulam stable. 
	\end{proof}
	
	\begin{theorem}
		Let $T$ be a densely defined closable operator  from 
		$D(T) \subset H$ to 
		$K$. Then  $C_{\overline T} = (I + T^{*} \overline T)^{-1}$ is Hyers-Ulam stable iff $T$ is  bounded.
	\end{theorem}
	\begin{proof}
		It is observed from \cite{MR2953553} that $C_{\overline T}$ is a positive self-adjoint and bounded operator in domain $H$ with $0 \leq C_{\overline T} \leq I$. So, $C_{\overline T}$ is closed and $R(C_{\overline T}) = D(T^{*} \overline T)$ is a core for $\overline T$. Suppose that  $C_{\overline T}$ is Hyers-Ulam stable. For given $x\in D(\overline T)$, there exists a sequence $\{x_{n}\}$ in $D(T^{*} \overline T)$ such that $x_{n} \to x$ and $\overline Tx_{n} \to \overline Tx$ as $n\to \infty$.  As $C_{\overline T}$ is Hyers-Ulam stable, $R(C_{\overline T})$ is closed. Hence $x\in D(T^{*} \overline T)$, thus  $D(\overline T) = D(T^{*} \overline T) = R(C_{\overline T})$. As $\overline T$ is closed in the closed domain $D(\overline T)$, $ T$ is bounded.
		
		Conversely, suppose that $T$ is both closable and bounded. Then, $\overline T$ is closed and bounded. Also,  $D(\overline T)$ is closed. $D(\overline T) = H$ because of $\overline{D(\overline T)} = H$. We have
		\begin{equation}\label{eqn_5}
			H\supset \overline { D(T^{*}\overline T)} \supset D(\overline T) = H.
		\end{equation}
		$T^{*} \overline T$ is closed because of the self-adjointness property of $T^{*} \overline T$ and $T^{*} \overline T$ is bounded because $T$ is bounded. So, $D(T^{*} \overline T)$ is closed. From (\ref{eqn_5}), we have $D(T^{*}\overline T) = H= R(C_{\overline T})$, hence $C_{\overline T}$ is Hyers-Ulam stable.
	\end{proof}
	\begin{theorem}\label{thm2.15}
		Let $T$ be a densely defined closable operator from $D(T)\subset H$ to $K$. If $T$ is Hyers-Ulam stable, then $Z_{\overline T} = \overline T{(I+ T^{*} \overline T)^{{-}\frac{1}{2}}} = \overline T{ {C_{\overline T}}^{\frac{1}{2}}}$ is Hyers-Ulam stable. 
	\end{theorem}
	\begin{proof}
		$Z_{\overline T}$ is bounded in domain $H$ with $\|Z_{\overline T}\| \leq 1$  implies that $Z_{\overline T}$ is closed . We claim that $R(Z_{\overline T}) = R(\overline T)$. By the definition of $Z_{\overline T}$, $R(Z_{\overline T}) \subset R(\overline T)$.
		\begin{equation}\label{equ 6}
 Z_{\overline T} (I + T^{*} \overline T)^{\frac{1}{2}}x= \overline Tx, \quad  \text{ for all } x\in D(I + T^{*} \overline T)^{\frac{1}{2}}.
		\end{equation}
		The equality (\ref{equ 6}) says that $D(Z_{\overline T} (I + T^{*}\overline T)^{\frac{1}{2} })\subset D(\overline T)$. Hence  $(Z_{\overline T} (I + T^{*} \overline T)^{-\frac{1}{2}})^{*} = (I+ T^{*} \overline T)^{-\frac{1}{2}} {Z_{\overline T}}^{*} $ because of the boundedness of $Z_{\overline T}$ and the self-adjointness property of the operator $(I + T^{*} \overline T)^{-\frac{1}{2}}$. 
		Let  $x\in D(\overline T)$ and  $y\in H$. Then $ (I+ T^{*} \overline T)^{-1}y \in D(T^{*} \overline T)$. Also, we have 
		\begin{equation}\label{equ 7}
			\begin{split}
				\langle ((I+ T^{*} \overline T)^{-\frac{1}{2}} {Z_{\overline T}}^{*} \overline T + (I+T^{*}\overline T)^{-1})x, \; y\rangle  
				&= \langle \overline Tx,\;  Z_{\overline T} (I+ T^{*} \overline T)^{-\frac{1}{2}}y\rangle + \langle x, \; (I + T^{*} \overline T)^{-1} y\rangle\\
				&=  \langle \overline Tx,\;  \overline T (I+ T^{*} \overline T)^{-1}y\rangle + \langle x, \; (I + T^{*} \overline T)^{-1} y\rangle\\
				&= \langle x,\;  T^{*} \overline T (I+ T^{*} \overline T)^{-1}y\rangle + \langle x, \; (I + T^{*} \overline T)^{-1} y\rangle\\
				&= \langle x, \; y\rangle.
			\end{split}
		\end{equation}
		
		From equality (\ref{equ 7}), we conclude  that $x= ({C_{\overline T}^{\frac{1}{2}}} {Z_{\overline T}}^{*} \overline T + C_{\overline T}) x= {C_{\overline T}}^{\frac{1}{2}} ({Z_{\overline T}}^{*} \overline T + C_{\overline T}^{\frac{1}{2}}) x$, which implies that $x\in R({C_{\overline T}}^{\frac{1}{2}})$. Hence, $D(\overline T)\subset  R({C_{\overline T}}^{\frac{1}{2}}) \subset D(\overline T)$, which confirms that $D(\overline T) =  R({C_{\overline T}}^{\frac{1}{2}})$. For arbitrary $u\in D(\overline T)$, there is an element $v\in H$ such that $u= {C_{\overline T}^{\frac{1}{2}}}v$, which shows that  $\overline Tu= \overline T{C_{\overline T}^{\frac{1}{2}}}v = Z_{\overline T}v$. Thus, $R(\overline T) \subset R(Z_{\overline T})$ implies $R(\overline T)= R(Z_{\overline T})$.
		Hence $Z_{\overline T}$ is Hyers-Ulam stable iff $\overline T$ is Hyers-Ulam stable. By Proposition \ref{pro 2},  $Z_{\overline T}$ is Hyers-Ulam stable. 
	\end{proof}

	\begin{theorem}\label{thm 5}
		Let $T$ be a densely defined closable operator from $D(T) \subset H$ to $K$. If $T$ is  Hyers-Ulam stable, then $\sigma({T^{*}\overline T}_{{N(T^{*}\overline T)}^{\perp}}) \subset [\,r, \infty)$, for some $r>0$. 		
		
	\end{theorem}
	\begin{proof}
		
		If  $N(\overline T)={0}$,  then  $\sigma (T^{*}\overline T) = \sigma(T^{*}\overline T_{{N(T^{*}\overline T)}^{\perp}})= \sigma(T^{*}\overline T_{N(\overline T)})\cup \sigma(T^{*}\overline T_{{N(T^{*}\overline T)}^{\perp}})$.  Hence
	\begin{equation}\label{equ 8}
	\sigma (T^{*}\overline T) =  \sigma(T^{*}\overline T_{N(\overline T)})\cup \sigma(T^{*}\overline T_{{N(T^{*}\overline T)}^{\perp}}).
	\end{equation}

		 \noindent If $N(\overline T)\neq {0}$, 		 
      we claim that $\lambda \notin \sigma({T^{*}\overline T}_{N({T^{*}\overline T})}) \cup \sigma(T^{*}\overline T_{{N(T^{*}\overline T)}^{\perp}})$, when $\lambda \notin \sigma(T^{*}\overline T)$. $(T^{*}\overline T -\lambda)$ is bijective and $(T^{*}\overline T -\lambda)^{-1}$ is bounded on domain $H$. Hence both $(T^{*}\overline T_{N(T^{*}\overline T)} - \lambda)$ and $(T^{*}\overline T_{{{N(T^{*}\overline T)}^{\perp}}} - \lambda)$ are injective.  Then, by Proposition \ref{pro 2},  $R(\overline T)$ is closed.  So, $R(T^{*}\overline T)$ is also closed and $R(T^{*}\overline T) = {N(T^{*}\overline T)}^{\perp}$.
		
		For any arbitrary $y\in N(T^{*}\overline T) = N(\overline T)$, there exists $x = x_{1} + x_{2} \in D(T^{*}\overline T)$, where $x_{1}\in N(T^{*}\overline T), x_{2}\in C(T^{*}\overline T)$ such that $(T^{*}\overline T - \lambda)(x_{1}+ x_{2}) =y$. So, 		
		$T^{*}\overline T x_{2} -\lambda x_{2} = \lambda x_{1} +y \in N(T^{*}\overline T) \cap {N(T^{*}\overline T)}^{\perp}$. Then $(T^{*}\overline T_{N(T^{*}\overline T)} - \lambda)x_{1} =y$, hence  $(T^{*}\overline T _{N(T^{*}\overline T)}- \lambda)$ is surjective.
		Similarly, for an arbitrary $z\in {N(T^{*}\overline T)}^{\perp}$, there exists $w= w_{1}+ w_{2} \in D(T^{*}\overline T)$, where $w_{1} \in N(T^{*}\overline T), w_{2}\in {N(T^{*}\overline T)}^{\perp} \cap D(T^{*}\overline T)$ such that $(T^{*}\overline T - \lambda)w =z$. Now we can write that $T^{*}\overline Tw_{2}-z-\lambda w_{2} = \lambda w_{1} \in {N(T^{*}\overline T)}^{\perp} \cap N(T^{*}\overline T)$. It shows that $T^{*}\overline T w_{2} -\lambda w_{2}= z$ implies $(T^{*}\overline T_{{{N(T^{*}\overline T)}^{\perp}}} - \lambda)$ is surjective. Thus $(T^{*}\overline T_{N(T^{*}\overline T)} -\lambda)^{-1}$ and $(T^{*}\overline T_{{{N(T^{*}\overline T)}^{\perp}}} -\lambda)^{-1}$ both are bounded. Hence, $\lambda \in \rho(T^{*}\overline T_{N(T^{*}\overline T)}) \cap \rho(T^{*}\overline T_{{N(T^{*}\overline T)}^{\perp}})$ implies $\lambda \notin \sigma (T^{*}\overline T_{N(T^{*}\overline T)}) \cup \sigma(T^{*}\overline T_{{N(T^{*}\overline T)}^{\perp}})$. It says that
		\begin{equation}\label{equ 9}
			\sigma (T^{*}\overline T_{N(T^{*}\overline T)}) \cup \sigma(T^{*}\overline T_{{N(T^{*}\overline T)}^{\perp}}) \subset \sigma(T^{*}\overline T).
		\end{equation}

		   Next, we claim that $T^{*}\overline T_{{N(T^{*}\overline T)}^{\perp}}$ is self-adjoint. Let $T^{*}\overline T =A$ and ${{N(T^{*} T)}^{\perp}}= N$. Then  $A_{N} =( A^{*})_{N} \subset A^{*} = A$, so $ A_{N} \subset (A_{N})^{*}$, hence $A_{N}$ is symmetric. Since $A$ is positive self-adjoint, there is a $\lambda \in \mathbb C \setminus \mathbb R$ with $R(A- \lambda) = H,  R(A- \overline \lambda) =H$ and $\lambda, \overline \lambda \in \rho(A)$. From (\ref{equ 8}) and (\ref{equ 9}), $\lambda, \overline \lambda  \notin \sigma(A)$ imply $ \lambda, \overline \lambda  \notin \sigma(A_{N})$. Therefore $R(A_{N} -\lambda) = N(T^{*}\overline T)^{\perp} = R(A_{N} - \overline \lambda)$. It has been proven in \cite{MR2953553} that any symmetric operator $S$ if there exists a complex number $\lambda$ such that $R(S-\lambda) =H$ and $R(S- \overline \lambda)$ is dense in $H$, then $S$ is self-adjoint. Thus, $A_{N}$ is also a self-adjoint and positive operator. The closedness of $R(T^{*}\overline T)$ guarantees that $0\notin \sigma(T^{*} \overline T_{{N(T^{*}\overline T)}^{\perp}}) = \sigma(A_{N})$. Therefore, there exists  $r>0$ such that $\sigma (T^{*}\overline T_{{N(T^{*}\overline T)}^{\perp}}) \subset [\, r, \infty)$.

	\end{proof}
 \begin{remark}
 The reverse inclusion of (\ref{equ 9}) is also true. Let $\lambda \in \sigma(T^{*}\overline T)$.  It implies that either $(T^{*}\overline T- \lambda) $ is not injective or not surjective.
		Suppose $(T^{*}\overline T - \lambda)$ is not injective. Then there exists $u\neq 0$ in $D(T^{*}\overline T)$ with $u= u_{1}+ u_{2}$, where $u_{1}\in N(T^{*}\overline T)$ and $u_{2} \in {N(T^{*}\overline T)}^{\perp}$ such that $(T^{*}\overline T- \lambda)u= 0$. So, $T^{*}\overline Tu_{2} -\lambda u_{2} =\lambda u_{1} \in N(T^{*}\overline T) \cap {N(T^{*}\overline T)}^{\perp},$ which guarantees that $u_{1}= 0$ and $(T^{*}\overline T_{{N(T^{*}\overline T)}^{\perp}} - \lambda)$ is not injective. Thus $\lambda \in \sigma(T^{*}\overline T_{{N(T^{*}\overline T)}^{\perp}})$.  
		
Suppose that  $(T^{*}\overline T -\lambda)$ is not surjective. Then there is an element $v\in H$ which has no pre-image with respect to $(T^{*}\overline T - \lambda)$. Let
$v= v_{1}+ v_{2}$, where $v_{1}\in N(T^{*}\overline T)$, $v_{2}\in {N(T^{*}\overline T)}^{\perp}$.  If both $(T^{*} \overline T_{N(T^{*}\overline T)} - \lambda)$ and $(T^{*}\overline T_{{N(T^{*}\overline T)}^{\perp}} -\lambda)$ are surjective, there exist $s_{1}\in N(T^{*}\overline T)$ and $s_{2}\in {N(T^{*}\overline T)}^{\perp} \cap D(T^{*}\overline T)$ such that $T^{*}\overline T_{N(T^{*}\overline T)} (s_{1}-\lambda s_{1} ) =v_{1}$ and $T^{*}\overline T_{{N(T^{*}\overline T)}^{\perp}}( s_{2}-\lambda s_{2}) =v_{2}$ which confirm the  pre-image of $v$. Therefore, both the operators $(T^{*} \overline T_{N(T^{*}\overline T)} - \lambda)$ and $(T^{*}\overline T_{{N(T^{*}\overline T)}^{\perp}} -\lambda)$ cannot be surjective at the same time. Therefore we have
		\begin{equation}\label{equ 10}
			\sigma (T^{*}\overline T) \subset \sigma(T^{*}\overline T_{N(T^{*}\overline T)})\cup \sigma(T^{*}\overline T_{{N(T^{*}\overline T)}^{\perp}}).
		\end{equation}
\end{remark}
	\begin{remark}
	\begin{enumerate}
		\item 			The converse of Theorem \ref{thm2.15} is  true when $N(\overline T) = \overline{N(T)}$. It can be proven by Corollary \ref{cor 1}.
	
	\item 	The converse of Theorem \ref{thm 5} is  true when $N(\overline T) = \overline{N(T)}$. 
	Indeed,  suppose that $\sigma({T^{*}\overline T}_{{N(T^{*}\overline T)}^{\perp}}) \subset [\,r, \infty)$ for some $r>0$ which implies $0 \in \rho(T^{*} \overline T_{{N(T^{*}\overline T)}^{\perp}})$. Then $(T^{*} \overline T_{{N(T^{*}\overline T)}^{\perp}})^{-1}$ exists and bounded with $D(T^{*} \overline T_{{N(T^{*}\overline T)}^{\perp}})^{-1} = {{N(T^{*} \overline T)}^{\perp}}$. Further,
	\begin{equation}
	R(T^{*} \overline T) = R(T^{*}\overline T_{N(T^{*}\overline T)}) +R(T^{*}\overline T_{{N(T^{*}\overline T)}^{\perp}}) =R(T^{*}\overline T_{{N(T^{*}\overline T)}^{\perp}}) = {{N(T^{*}\overline T)}^{\perp}}.
	\end{equation}
	Thus $R(T^{*}\overline T)$ is closed which implies $R(\overline T)$ is closed \cite{MR2914387}.  Therefore, by Corollary  \ref{cor 1},  $T$ is Hyers-Ulam stable.

\item  It is known that  non-zero complex number  $\lambda$ lies in  $ \rho(T^{*}\overline T _{N(T^{*} \overline T)})$.		If a densely defined closable operator $T$ from $D(T) \subset H$ to $K$ is Hyers-Ulam stable, then $\sigma(T^{*} \overline T) \subset \{0\} \cup [\, r, \infty)$.

	\end{enumerate}
	\end{remark}
	\begin{example}
		
		Define $T$ on $\ell^{2}$ by $T(x) =( 0, 2x_{2}, 3x_{3},\ldots,nx_{n},\ldots)$, where $x=(x_{1}, x_{2}, x_{3},\ldots,x_{n},\ldots)\in \ell^{2}$ and $D(T)=\{x\in \ell^{2} : ( 0, 2x_{2}, 3x_{3},\ldots,nx_{n},\ldots )\in \ell^{2}\}$. Then, $T$ is a self-adjoint operator. Moreover,  $\sigma(T)= \{(n-1): n\in \mathbb N \setminus \{2\}\}$ and $\sigma(T^{2}) = \sigma (T^{*}T)= \{(n-1)^{2}: n\in \mathbb N\setminus \{2\}\}$. Then $\sigma (T^{*}T_{{N(T^{*}T)}^{\perp}}) = \{n^{2} : n\in \mathbb N \setminus \{1\}\} \subset [4, \infty)$. Therefore, $T$ is Hyers-Ulam stable. 
		
	\end{example}
	\begin{example}
		Define $T: D(T) \subset \ell^{2} \to \ell^{2}$ by $T(x)= (x_{1}, 2x_{2}, \frac{x_{3}}{3}, 4x_{4}, \frac{x_{5}}{5},\ldots)$, where $D(T)= \{(x_{1}, x_{2}, x_{3},\ldots,x_{n},\ldots)  : (x_{1}, 2x_{2}, \frac{x_{3}}{3}, 4x_{4}, \frac{x_{5}}{5},\ldots)\in \ell^{2} \}$. 
		Then $T$ is closed and $N(T^{*})= N(T^{*}T) =N(T) =\{0\}$. So, $\sigma(T^{*}T)=\sigma(T^{*}T_{{N(T^{*}T)}^{\perp}})$. It is easy to observe that $\{\frac{1}{(2n +1)^{2}} : n\in \mathbb N\} \subset \sigma_{p}(T^{*}T)$. Thus $0\in \sigma(T^{*}T) = \sigma(T^{*}T_{{N(T^{*}T)}^{\perp}})$. Therefore, $T$ is not Hyers-Ulam stable.
	\end{example}
	\begin{corollary}
		Let $T$ be a densely defined closable Hyers-Ulam stable operator from $D(T) \subset H$ to $H$. If $ T$ is essentially self-adjoint and  $N(\overline T) =\{0\}$,  then $T^{2}$ is Hyers-Ulam stable.  
	\end{corollary}
	\begin{proof}
		Since $\overline T$ is the self-adjoint operator,  $N(T^{*}\overline T) = N({\overline T}^{2}) = N(\overline T) =\{0\}$. By Theorem \ref{thm 5},  $\sigma(T^{*}\overline T_{N(\overline T)}) =\emptyset$ and $\sigma(T^{*}\overline T) = \sigma(T^{*}\overline T_{{N(T^{*}\overline T)}^{\perp}}) \subset [\, r, \infty)$, for some $r>0$. Hence $0\in \rho (T^{*}\overline T) =\rho({\overline T}^{2})$ implies $R({\overline T}^{2}) = H$. Also, $\overline T^{2}= T^{*} \overline T$ is self-adjoint. Thus, $\overline T^{2}$ is also Hyers-Ulam stable.  Moreover, $N(T^{2}) = N( T^{*} \overline T) = \{0\}$  with $T^{2} \subset T^{*}\overline T$. From the definition of the Hyers-Ulam stability of $T^{*}\overline T$, we can conclude that $T^{2}$ is Hyers-Ulam stable.
	\end{proof}
	\begin{proposition}
		Let $S$ and $T$ be two densely defined Hyers-Ulam stable operators from $D(S) \subset H$ to $K$ and $D(T) \subset K$ to $H$ respectively with $T \subset S^{*}$ and $S \subset T^{*}$. Assume that  $TS$ is a densely defined operator and dim$ (N(\overline T)) < \infty $. If $\overline {TS}$ is self-adjoint and $N(TS)^{*} = \overline {N(TS)}$, then ${TS}$ is Hyers-Ulam stable.
	\end{proposition}
	\begin{proof}  The operators 
		$\overline T$ and $\overline S$ both are Hyers-Ulam stable because of the Hyers-Ulam stability of $T$ and $S$. Also,  $R(\overline T)$ and $R(\overline S)$ both are closed. Then ${\bar T}  {\bar S}$ is closed, and $R(\bar T \bar S)$ is also closed because $\bar T$ and $\bar S$ both are normally solvable with the given condition dim$N(\bar T) <\infty$ implies ${\bar{T}  \bar {S}}$ is also normally solvable.
		\begin{equation}\label{equ 12}
			\overline {TS} \subset \bar T  \bar S \subset S^{*}T^{*} \subset (\bar T  \bar S)^{*} \implies \overline {TS} \subset \bar T  \bar S \subset (\overline {TS})^{*}.
		\end{equation}
		Since $\overline {TS}$ is self-adjoint,  $\overline {TS} = \bar T  \bar S$. Then, the Hyers-Ulam stability of $\overline {TS}$  with the given condition $N(TS)^{*} = \overline {N(TS)}$ confirms the Hyers-Ulam stability of $TS$.
		
	\end{proof}
	\begin{proposition}
		Let $T$ be a positive self-adjoint Hyers-Ulam stable operator from $D(T) \subset H$ to $H$. Then $T^{\frac{1}{2^{n}}}$ is also Hyers-Ulam stable, for all $n\in \mathbb N$. Conversely,  if $T^{\frac{1}{2^{n}}}$ is Hyers-Ulam stable, for some $n\in \mathbb N$, then 		
		 $T$ is Hyers-Ulam stable. 
	\end{proposition}
	\begin{proof}
		Since $T$ is a self-adjoint operator,  $T$ is closed. Also $(T^{\frac{1}{2}})^{2}= T$, where $T^{\frac{1}{2}}$ is an unique square root of $T$ and  positive self-adjoint operator. It is proven in\cite{MR2914387} that $R(T)$ is closed iff $R(T) = R(T^{\frac{1}{2}})$, when $T$ is a positive self-adjoint operator. Hence, $T^{\frac{1}{2}}$ is also Hyers-Ulam stable. By the induction hypothesis, we can say that $T^{\frac{1}{2^{n}}}$ is also Hyers-Ulam stable for all $n\in \mathbb N$.
    
    Conversely, suppose that  $T^{\frac{1}{2^{n}}}$ is Hyers-Ulam stable, for some  $n\in \mathbb N$. We claim that $T^{\frac{1}{2^{n-1}}}$ is also Hyers-Ulam stable. We know that $ D(T^{\frac{1}{2^{n-1}}}) \subset D(T^{\frac{1}{2^{n}}})$ and $N(T^{\frac{1}{2^{n}}}) = N(T^{\frac{1}{2^{n-1}}})$ because $T^{\frac{1}{2^{n-1}}}$ is a self-adjont operator. Now we consider an arbitrary $\varepsilon \geq 0$ with $x \in D(T^{\frac{1}{2^{n-1}}})$ such that $\|T^{\frac{1}{2^{n-1}}}(x)\| = \|(T^{\frac{1}{2^{n}}})^{2}(x)\| \leq \varepsilon$. Since $T^{\frac{1}{2^{n}}}$ is Hyers-Ulam stable, so we get an element $x_{1}\in N(T^{\frac{1}{2^{n}}})$ such that $\|(T^{\frac{1}{2^{n}}})(x) - x_{1}\| \leq M_{T^{\frac{1}{2^{n}}}} \varepsilon $. Now this inequality $\|(T^{\frac{1}{2^{n}}})(x) - x_{1}\|^{2} \leq (M_{T^{\frac{1}{2^{n}}}} \varepsilon)^{2} $ confirms that the inequality $\|(T^{\frac{1}{2^{n}}})(x)\| \leq M_{T^{\frac{1}{2^{n}}}} \varepsilon $. Again we get an element $x_{0}\in  N(T^{\frac{1}{2^{n}}})$ such that $\|x- x_{0}\| \leq (M_{T^{\frac{1}{2^{n}}}})^{2} \varepsilon $.  Thus  $T^{\frac{1}{2^{n-1}}}$ is Hyers-Ulam stable. By the induction hypothesis, it can be shown that $T$ is Hyers-Ulam stable.
	\end{proof}
 \begin{definition}
	Let $T$ be an operator on $H$.  The operator $T$ is said to be a paranormal  operator if $\|Tx\|^{2} \leq \|T^{2}x\| \|x\|,$ for all $x\in D(T^{2})$.
\end{definition}
	\begin{theorem}\label{thm 4}
		Let $T$ be a paranormal operator on $H$. If $0\in \rho(T)$,  then $T^{n}$ is always Hyers-Ulam stable, for each  $n\in \mathbb N$.
	\end{theorem}
	\begin{proof}
	   As $\sigma (T) \neq \mathbb C$, $T$ is closed.  It is proved in \cite{MR4557345} that $N(T) = N(T^{n})$,  for all $n\in \mathbb N$. So $N(T^{n}) = 0$ which confirms the existence of $(T^{n})^{-1}$. Hence  $(T^{-1})^{n}$ is bounded with $(T^{n}) (T^{-1})^{n} = I_{H}$ and $(T^{-1})^{n} (T^{n})= I_{D(T^{n})} \subset I_{H}$. The left and right inverses of $T^{n}$ exist. Thus, the inverse of $T^{n}$ exists and is unique \cite{MR3770359}.  Also $(T^{n})^{-1} = (T^{-1})^{n}$. Then $R(T^{n}) =H $ and $0\in \rho(T^{n})$.  Moreover,  $T^{n}$ is closed. Therefore  $T^{n}$ is Hyers-Ulam stable for all $n\in \mathbb N$.   
	\end{proof}
 \begin{remark}
 Theorem \ref{thm 4} holds good for
 self-adjoint, quasi-normal, subnormal and hyponormal operators because they are paranormal operators.
 \end{remark}
 \begin{corollary}
	Let $T$ be a symmetric operator on $H$. If $0\in \rho(T)$, then $T^{n}$ is Hyers-Ulam stable, for all $n\in \mathbb N$.
	\end{corollary}
	\begin{proof}Since $T \subset T^{*}$ and $R(T) = H$, $T$ is self-adjoint. From  Theorem \ref{thm 4}, we conclude that $T^{n}$ is Hyers-Ulam stable.
	\end{proof}

\section{Sum and Product of Hyers-Ulam Stable Operators}
 In this section, we establish various results to show that the sum and product of two Hyers-Ulam stable to be Hyers-Ulam stable and find some necessary and sufficient conditions of $R(\overline{\mathcal A})$ to be closed and $\overline{\mathcal A}$ is Hyers-Ulam stable, where $\overline{\mathcal A}$ is the closure of the $2 \times 2$ block matrix $\mathcal A$.

\color{black}
\begin{theorem}\label{thm 7}
	Let $S$ be a closable operator from $D(S) \subset H$ to $K$. Let $T$ be a  closable operator with the Hyers-Ulam stability from $D(T) \subset H$ to $K$. If $D(T) \subset D(S)$ and $\|Sx\| \leq b \|Tx\|$ for all $x\in D(T)$, where $0 \leq b< 1 $, then $\overline {(S+T)}$ is  Hyers-Ulam stable. Moreover, $M_{\overline {S+T}} \leq \frac{1}{1-b} M_{\overline T}$.
\end{theorem}
\begin{proof}
	It is shown in \cite{MR2857980} that $\| {\overline S}x\| \leq b \|{\overline T} x\|$, for all $x\in D(\overline T) \subset D(\overline S)$, and $\overline {(S + T)} = \overline S + \overline T$, hence $(\overline S+\overline T)$ is closed.  We will show that $R(\overline S+\overline T)$ is closed in $K$.
	We know that $R(\overline S+\overline T)$ is closed iff $\gamma(\overline S+ \overline T) > 0 $. So, it suffices to show that $\gamma(\overline S+\overline T) >0$.  Note that 
	$N(\overline T)\subset N(\overline S) $ and $N(\overline T) = N(\overline T) \cap N(\overline S) \subset N(\overline T+\overline S)$.
	Thus, $C(\overline T+\overline S) ={N(\overline T+\overline S)}^{\perp} \cap D(\overline T+\overline S) = {N(\overline T+\overline S)^{\perp}} \cap D(\overline T) \subset {N(\overline T)}^{\perp} \cap D(\overline T)= C(\overline T)$.
	
	\noindent Moreover, $dist(x, N(\overline T)) \geq dist (x, N(\overline T+\overline S))$, for all $x\in {N(\overline T+\overline S)^{\perp}} \cap D(\overline T+\overline S)$.
	Then we get the following inequality:
	\begin{align*}
		\frac{(1-b) \|\overline Tx\|}{dist(x, N(\overline T))} \leq \frac{\|(\overline T+\overline S)x\|}{dist (x, N(\overline T))} \leq \frac{\|(\overline T+\overline S) x\|}{dist(x, N(\overline T+\overline S))}, \text{ for all } x\in C(\overline T+\overline S).
	\end{align*}
	From the above inequality, we have 
	\begin{align*}
		(1-b)  {\inf_{x\in C(\overline T)} }\frac{ \|\overline Tx\|}{dist (x, N(\overline T))} \leq \inf_{x\in C(\overline T+\overline S) } {\frac{\|(\overline T+\overline S)x\|}{dist(x, N(\overline T+\overline S))}}.
	\end{align*}
	Hence, $(1-b)\gamma (\overline T) \leq \gamma(\overline T+\overline S) $. As $\gamma(\overline T)>0$ and $0 \leq b<1$, $\gamma(\overline T+\overline S) >0$.
	Therefore, $R(\overline T+\overline S)$ is closed by Theorem \ref{thm 1}  and $(\overline T+\overline S) = \overline {T+S}$ is Hyers-Ulam stable with 
	\begin{align*}
		\frac{1}{\gamma(\overline {T+S})} \leq \frac{1}{(1-b) \gamma(\overline T)} \text{ \  implies  \ }  M_{(\overline {S+T})} \leq \frac{1}{(1-b)} M_{\overline T}.
	\end{align*}
	
\end{proof}
Now, we will prove the above Theorem \ref{thm 7} without assuming the closableness of $S$ and $T$. 
\begin{theorem}\label{thm 8}
    Let $S$ and $T$ be two operators from $H$ into $K$ with $D(T) \subset D(S)$ and $\|Sx\| \leq b \| Tx\|$, for all $x\in D(T), 0 \leq b < 1$. If $T$ is Hyers-Ulam stable, then $S+T$ is Hyers-Ulam stable with $M_{S+T} \leq \frac{M_{T}}{(1-b)}$.
\end{theorem}
\begin{proof}
 We have $N(T) \subset N(T+S)$ because  $N(T) \subset N(S)$. For given $\varepsilon > 0$ and $x\in D(S+T)$ with $\|(S +T)x\| \leq \varepsilon $, we get 
\begin{equation}\label{equ 13}
    (1-b)\|Tx\| \leq \|Tx\| -\|Sx\|\leq \|(S+T)x\|\leq \varepsilon .
\end{equation}
Then $\|Tx\| \leq \frac{\varepsilon}{1-b}$. Since $T$ is Hyers-Ulam stable, then there exists $M \geq 0$ with  $x_{0} \in N(T)$ such that $\|x-x_{0}\| \leq \frac{M \varepsilon}{1-b}$.
 Therefore, $S+T$ is Hyers-Ulam stable with $M_{S+T} \leq \frac{M_{T}}{(1-b)}$.
\end{proof}
\begin{theorem} Let $S$ be a closable operator from $D(S) \subset H$ to $K$. Let $T$ be a closable operator with the Hyers-Ulam stability from $D(T) \subset H$ to $K$. Assume $D(T) \subset D(S)$ and $\|Sx\| \leq b \|Tx\|$ for all $x\in D(T)$, where $0 \leq b< 1 $. If $N(\overline T)= \overline{N(T)}$, then $$M_{S+T} = M_{\overline {S+T}} \leq \frac{1}{1-b} M_{\overline T} = \frac{1}{1-b} M_{ T}.$$
\end{theorem}
\begin{proof} 
It has been proven in \cite{MR2857980} that $\overline {S+T} = \overline S + \overline T$. By Theorem \ref{thm 8}, $S+T$ is Hyers-Ulam stable. Hence, both $\overline {S+T}$ and $\overline T$  are also Hyers-Ulam stable. Note that $M_{\overline T} = M_{T}$ because $N(\overline T)= \overline{N(T)}$. We can easily show that $N{\overline {(S+T)}} = \overline {N(S + T)}$ from the given condition $N(\overline T)= \overline{N(T)}$. By Corollary \ref{cor 1}, we get that  $M_{S+T} = M_{\overline{S+T}}$. Therefore by Theorem \ref{thm 8}, we have the following inequality $$M_{S+T} = M_{\overline {S+T}} \leq \frac{1}{1-b} M_{\overline T} = \frac{1}{1-b} M_{ T}.$$
\end{proof}

\begin{theorem}
	Let $T$ be a closable operator from $D(T) \subset H$ to $K$ with an assumption that $D(T)$ is the decomposable domain. Let $T$ and $T^{\dagger}$ be Hyers-Ulam stable operators. Assume that the closable operator $S$ is a $T$-bounded operator with $D(S) =D(T)$ and $\|Sx\| \leq b\|Tx\|$, for some non-negative $b\in \mathbb R$. If $\| S\| \|\overline T^{\dagger}\| < 1$, then $\overline {(S+T)} = \overline S + \overline T$ is also Hyers-Ulam stable. Moreover, $S + T$ is Hyers-Ulam stable when $N(\overline{S + T}) = \overline{N(S + T)}$.
\end{theorem}
\begin{proof}
Theorem \ref{thm 6} clearly states that $\overline T$ is bounded.  Then $D(\overline T)$ is closed and $\|\overline Sx\| \leq b \|\overline Tx\|$, for all $x\in D(\overline T) \subset D(\overline S)$. Also, $\overline S$ is also bounded on domain $D(\overline T)$. Therefore $(\overline S+\overline T)$ is bounded. From \cite{MR2857980}, we know that $\overline S + \overline T \subset \overline {S+T}$. Note that  $\overline S + \overline T$ is closed because $D(\overline T)$ is closed, and $\overline S + \overline T$ is bounded. Thus $\overline S + \overline T = \overline {S+T}$. 
	
	Moreover, $R(\overline T)$ is closed because of the Hyers-Ulam stability of $T$. So, $\overline T^{\dagger}$ is bounded. Then $\gamma(\overline T) = \frac{1}{\|\overline T^{\dagger}\|} > 0$. Hence $   \|\overline Tx\| \geq \gamma(\overline T) \|x\|$  for all $x\in C(\overline T)$. Therefore for each $x\in C(\overline T)$, we have
\begin{equation}\label{equ 14}
		\|(\overline S+\overline T)x\| \geq \gamma(\overline T) \|x\| - \|S\| \|x\| = \frac{1}{\|\overline T^{\dagger}\|} \|x\| - \|S\| \|x\| =\frac{(1- \|\overline T^{\dagger}\| \|S\|)}{\| \overline T^{\dagger}\|} \|x\|.
\end{equation}

\begin{equation}\label{equ 15}
		N(\overline T) \subset N(\overline S) \text{ and }   N(\overline T+ \overline S) \supset  N(\overline T)   \text { \ implies \ } N(\overline T)^{\perp}  \supset N(\overline T+\overline S)^{\perp}. 
\end{equation}
	The above two inequalities (\ref{equ 14}) and (\ref{equ 15}) say that	$$\| (\overline S+ \overline T)x\| \geq \frac{(1- \|\overline T^{\dagger}\| \|S\|)}{\| \overline T^{\dagger}\|} \|x\|,  \text { for all } x\in C(\overline S+\overline T).$$ 
	Since $\|S\| \|\overline T^{\dagger}\| < 1$, 
	$\frac{(1- \|\overline T^{\dagger}\| \|S\|)}{\|\overline T^{\dagger}\|} >0$. 	Hence  $\gamma{{(\overline S+\overline T)} }> 0$  and by Theorem \ref{thm 1},  $R(\overline S+\overline T)$ is closed, so  $\overline S + \overline T$ is Hyers-Ulam stable. Furthermore, $S + T$ is Hyers-Ulam stable by Corollary \ref{cor 1}.
\end{proof}
\begin{theorem}
	Let $T$ be a closable Hyers-Ulam stable operator from $H$ to $K$ with the decomposable domain $D(T)$. Assume that $S$ is a $T$-bounded Hyers-Ulam stable operator with $D(S)= D(T)$ and $\|Sx\| \leq b\|Tx\|$, for all $x\in D(S)$ and  some $b\geq 0$. If $T^{\dagger}$ is Hyers-Ulam stable and $R(S) \subset R(T)^{\perp}$, then $\overline {(S+T)} = \overline S + \overline T$ is Hyers-Ulam stable. Moreover, if $N(\overline{S + T}) = \overline{N(S+T)}$, then $S + T$ is Hyers-Ulam stable.
\end{theorem}
\begin{proof}
By Theorem \ref{thm 6},  $\overline T$ is bounded, and $\| {\overline S}x\| \leq b \| {\overline T}x\|$, for all $x\in D(\overline T)$. Then $(\overline S+\overline T)$ is closed because $(\overline S+\overline T)$ is bounded on the closed domain $D(\overline T)$. Since $R(S) \subset R(T)^{\perp}$, $R(\overline S) \subset R(T)^{\perp}$ and  $R(T)^{\perp} = R(\overline T)^{\perp}$.  From \cite{MR2857980}, we know that $\overline S + \overline T \subset \overline{S+T}$, and the closedness of $\overline S + \overline T$ gives that $\overline{S+T} \subset \overline S + \overline T $. Hence,  $\overline S + \overline T = \overline{S+T}$. 	The Hyers-Ulam stability of $\overline S$ and $\overline T$ tell that both  $R(\overline S)$ and $R(\overline T)$  are closed in $K$. 
	
	We now claim that $R(\overline S+\overline T)$ is closed. 	Let $y \in \overline{R(\overline S+\overline T)}$. Then there exists a sequence $\{x_{n} \}$ in $D(\overline T)$ such that $(\overline S+\overline T) x_{n} \to y$ as $n\to \infty$.  As $R(\overline S) \subset R(\overline T)^{\perp}$,    both $\{\overline S{x_{n}}\}$ and $\{\overline T{x_{n}}\} $  are Cauchy sequences. Now, we can assume that $\overline S{x_{n}} \to z_{s}$ and $\overline T{x_{n}} \to z_{t}$ as $n\to \infty$, for some $z_{s}\in R(\overline S)$ and $z_{t} \in R(\overline T)$. Then there exists an element, say $x\in D(\overline T)$, such that $z_{t} =\overline Tx$. We get $\overline S{x_{n}} \to \overline Sx$ from the given inequality. Hence, $\overline Sx = z_{s}$ and $y= z_{s}+z_{t} = (\overline S+\overline T)x \in R(\overline S+\overline T)$. Therefore, $R(\overline{S+T})$  is closed and $\overline S + \overline T = \overline{S+T}$ is Hyers-Ulam stable. Furthermore, $S + T$ is Hyers-Ulam stable by Corollary \ref{cor 1}.

\end{proof} 
\begin{theorem}
	Let $T$ and $S$ be operators from $H$ to $K$. If $S$ is a closable operator and $T$ is bounded with $a\|x\| + \|Sx\| \leq \|Tx\|$, for all $x\in D(S) \subset H= D(T)$, for some $a>0$, then ${\overline {T+ S}} = T + \overline S$ is Hyers-Ulam stable.
\end{theorem}
\begin{proof}
From the given inequality, $\overline S$ is bounded in $D(\overline S)$. Hence, $T+\overline S$ is bounded and closed on the domain $D(\overline S)$. So, $D(\overline S)$ is closed and $T + \overline S = \overline {T+S}$.  It is enough to show that $R(T+\overline S)$ is closed. 
	
Let $y\in \overline{R(T+\overline S) }$. Then there exists a sequence $\{y_{n}\}$ in $ R(T+\overline S)$ such that $y_{n} \to y$ as $n \to \infty$. Thus, there exists a sequence $\{x_n\}$ in $D(\overline S)$ with $ (T+\overline S) x_{n} = y_{n}$, for all $n\in \mathbb N$. The given inequality is written of the form $$a\|x\| \leq \|Tx\| -\|\overline Sx\| \leq \|(T+\overline S) x\|.$$ Hence $\{x_n\}$ is a Cauchy sequence in the Hilbert space $D(\overline S)$. Thus   ${x_n} \to x$ as $n\to \infty $, for some $x\in D(\overline S)$.  As 
	 $\{y_n\} = \{(T+\overline S) x_n\} \to (T+\overline S) x$ as $n\to \infty$, $y= (T+\overline S) x\in R(T+\overline S)$.  Therefore $\overline {R(T+\overline S)} = R(T+\overline S) $, which shows that $R(T+\overline S)$ is closed. Thus ${\overline {T+ S}} = T + \overline S$ is  Hyers-Ulam stable.
\end{proof}
Let $H$ and $K$ be two Hilbert spaces. The product space $H \times K$ is a vector space if the linear operation is defined by  $\alpha (h_{1}, k_{1}) + \beta (h_{2}, k_{2}) = (\alpha h_{1} + \beta h_{2}, \alpha k_{1} + \beta k_{2})$, for all $h_{1}, h_{2}\in H$,  $k_{1}, k_{2}\in K,  \alpha, \beta \in \mathbb C $. Moreover, it is a  Hilbert space with an inner product $$\langle (h_{1}, k_{1}), (h_{2}, k_{2})\rangle = \langle h_{1}, h_{2}\rangle + \langle k_{1}, k_{2} \rangle, \text{ for } h_{1}, h_{2}\in H, k_{1}, k_{2}\in K.$$  Let $T$ and $S$ be two operators on $H$ and $K$ respectively with $(T \times S) (h,k)= (Th, Sk)$, for all $(h,k) \in D(T \times S)= \{(h,k)\in H \times K: (Th, Sk)\in H \times K\} = D(T) \times D(S)$.  We now discuss the Hyers-Ulam stability of $T \times S$.

\begin{theorem}\label{thm 15}
	Let $T$ and $S$ be two closable Hyers-Ulam stable operators on $H$ and $K$, respectively. Then $\overline {T \times S} = \overline T \times \overline S$ is Hyers-Ulam stable. If  $N(\overline{T \times S}) = \overline {N(T \times S)}$, then $T \times S$ is Hyers-Ulam stable.
\end{theorem}
\begin{proof}
	   It is obvious that $T \times S \subset \overline T \times \overline S$. Now, we will show that $\overline T \times \overline S$ is closed. let, $((u,v), (x,y)) \in \overline{G(\overline T \times \overline S)}$. There exists a sequence $\{(u_{n},v_{n})\}$ in $D(\overline T) \times D(\overline S)$ such that $(u_{n}, v_{n}) \to (u,v)$ and $(\overline{T}u_{n}, \overline{S}v_{n}) \to (x,y)$ as $n \to \infty$. Then $u\in D(\overline T), v\in D(\overline S)$ and $\overline{T}u =x, \overline{S}v = y$. Hence $\overline{T} \times \overline{S}$ is closed. Therefore $T \times S \subset \overline{T \times S} \subset \overline{T} \times \overline{S}$, where $\overline{T \times S}$ is the closure of $T \times S$. It is easy to show that $D(\overline T \times \overline S) \subset D(\overline{T \times S})$. Thus, $\overline{T \times S} = \overline T \times \overline S$. We know that $R(\overline{T \times S}) \subset \overline{R(T \times S)}$.

	   Now we claim that $R(\overline{T \times S})$ is closed. Let $(h',k') \in \overline{R(\overline{T \times S})} \subset \overline{R(T \times S)}$. Then there exists a sequence $\{(h'_{n},k'_{n})\}$ in $D(T \times S)$ such that $(Th'_{n}, Sk'_{n}) \to (h', k')$ as $n\to \infty$. The Hyers-Ulam stability of $T$ and $S$ tell that  $\overline {R(T)} \subset R(\overline T)$ and $\overline {R(S)} \subset R(\overline S)$. Let $h'\in R(\overline T)$ and  $k' \in R(\overline S)$. There exist $h\in D(\overline{T})$ and $k \in D(\overline S)$ with $h' = \overline{T}h$, $k'= \overline{S} k$. Thus $(\overline {T \times S})(h,k)= (h',k') \in R(\overline {T \times S})$. So, $R(\overline{T \times S})$ is closed. Therefore, $\overline{T \times S} = \overline T \times \overline S$ is  Hyers-Ulam stable.
   By Corollary \ref{cor 1}, we get that $T \times S$ is Hyers-Ulam stable.
\end{proof}
\begin{remark}  We now give an alternate idea to prove Theorem \ref{thm 15} by considering 
a $2 \times 2$ block matrix 	   of the form  $\mathcal {A} = \begin{bmatrix}
 T & 0\\
 0 & S
     
 \end{bmatrix}$. 
 Then, it is shown in \cite{MR2463978} that $\mathcal A$ is closable because  $T$ and $S$ are diagonally dominant with $T$-bound and $S$-bound less than $1$. Moreover, $\overline{\mathcal A} = \begin{bmatrix}
     \overline{T} & 0\\
     0 & \overline{S}
 \end{bmatrix}$
 and $R(\overline {\mathcal A}) = R(\overline {T} \times \overline{S}) = R(\overline{T}) \times R(\overline {S})$. As $T$ and $S$ are Hyers-Ulam stable, $R(\overline{T}) \times R(\overline {S})$ is closed. Therefore, the Hyers-Ulam stability of $\overline{\mathcal A}$ says the Hyers-Ulam stability of $\overline{T \times S}= \overline T \times \overline S$. Thus by Corollary \ref{cor 1},  $T \times S$ is Hyers-Ulam stable.

Let $\mathcal{A} = \begin{bmatrix}
  A & B\\
  C & E
  \end{bmatrix}$.  Here $A$ and  $E$ operators on $H$ and $K$ respectively whereas $B$ and $C$ are  operators from $K$ to $H$ and $H$ to $K$ respectively.  Note that  $$D(\mathcal A) = (D(A)\cap D(C)) \oplus (D(B) \cap D(E)).$$   The operators $A, B, C, E$ are used for $\mathcal{A}$ as explained above in the sequel. 
  We now provide an example of a $2 \times 2$ block matrix  $\mathcal {A}$, which is closed, but $R(\mathcal A)$ is not closed. 
\end{remark}
  \begin{example}
     Let $\mathcal A$ be an operator in $L^{2} (0,1) \times L^{2} (0,1)$.  We take $A = E = 0$ and $B = C= T_{\phi}$, where $T_{\phi}(f) = \phi f$ in the domain $D(T_{\phi}) =\{ f \in L^{2}(0,1): \phi f \in L^{2}(0,1)\} $ and $\phi(x) = x$, for all $x\in (0,1)$.  Note that     
      $R(T_{\phi})$ is a proper dense subset of $L^2(0,1)$ but $R(A) = R(E)= \{0\}$ is closed in $L^{2}(0,1)$.  Hence $\mathcal A$ is closed but $R(\mathcal A) = R(T_{\phi}) \times R(T_{\phi}) $ is also a proper dense subspace of $L^{2}(0,1) \times L^{2}(0,1)$.  Thus, $R(\mathcal A)$ is not closed.
  \end{example}
  \begin{theorem}
      Let $\mathcal A = \begin{bmatrix}
      A & B\\
      C & E
      \end{bmatrix}$ be an operator on $H \times K$ with $$\| Cx\| \leq b_{C} \|Ax\|, \text{ for all } x\in D(A) \subset D(C) \text{  and } 0 \leq b_{C} < 1$$ and $$\|B x\| \leq b_{B} \| Ex\|, \text{  for all } x\in D(E) \subset D(B) \text{ and } 0 \leq b_{B} <1.$$ If $A$ and $E$ both are Hyers-Ulam stable closable operators on $H$ and $K$ respectively,       
       then $R(\overline {\mathcal{A}})$ is closed. The converse is true when  $N(\overline A) = \overline{N(A)}$ and $N(\overline E) = \overline{N(E)}$.
  \end{theorem}
  \begin{proof}
      As shown in \cite{MR2463978}, we have that $\mathcal A$ is closable in $D(\mathcal A)= D(A) \oplus D(E)$ with $\overline {\mathcal{A}} = \begin{bmatrix}
          \overline A & B\\
          C & \overline E
      \end{bmatrix}$. For given $\varepsilon >0$ and $(x,y)^{T} \in D(\mathcal A)$ with $\| \mathcal{A} (x,y)^{T}\| < \varepsilon$, we get $\|Ax\| < \frac{2\varepsilon}{1- b_{C}}$ and $\|Ey\| < \frac{2\varepsilon}{1-b_{B}}$. Since $A$ and $E$  are Hyers-Ulam stable,  there exist two positive numbers, say $M_{1}$ and $M_{2}$ with two elements $x_{0}\in N(A) \subset N( C)$ and $y_{0} \in N(E) \subset N(B) $ such that $$\|x - x_{0}\| \leq \frac{2\varepsilon M_{1} }{1-b_{C}} \quad \text{ and } \|y - y_{0}\| \leq \frac{2\varepsilon M_{2} }{1-b_{B}}.$$ Then $\|(x,y)^{T} - (x_{0}, y_{0})^{T}\| \leq \varepsilon (\frac{2M_{1}}{1-b_{C}} +\frac{2M_{2}}{1-b_{B}})$, which shows that the Hyers-Ulam stability of  $\mathcal A$. By Proposition \ref{pro 2}, $\overline{\mathcal A}$ is also Hyers-Ulam stable. Therefore, $R(\overline{\mathcal A})$ is closed.
      
To prove the converse part, we will show that $A$ and $E$ are Hyers-Ulam stable. For given $\varepsilon >0$, $x \in D(A), y\in D(E)$ with $\|Ax\| \leq \varepsilon, \|Ey\| \leq \varepsilon$, we get $$\|\overline {\mathcal{A}}(x,y)^{T}\| \leq \varepsilon\sqrt{((1+ b_{B})^{2} + (1 + b_{C})^{2})}.$$ As $\overline{\mathcal A}$ is Hyers-Ulam stable, there exists an element $(u,v)^{T}\in N(\overline {\mathcal A})$  such that $\|(x,y)^{T} - (u,v)^{T}\| < \varepsilon M \sqrt{((1+ b_{B})^{2} + (1 + b_{C})^{2})}$, for some $M > M_{\overline{\mathcal A}}$. From $(u,v)^{T}\in N(\overline {\mathcal A)}$, $N(\overline A) = \overline{N(A)}$ and $N(\overline E) = \overline{N(E)}$,  we get  say $u_{0}\in N(A)$ and $v_{0}\in N(E)$ such that $\|x - u_{0}\| \leq \varepsilon M \sqrt{((1+ b_{B})^{2} + (1 + b_{C})^{2})}$ and $\|y - v_{0}\| \leq \varepsilon M \sqrt{((1+ b_{B})^{2} + (1 + b_{C})^{2})}$. Therefore, $A$ and $E$ are Hyers-Ulam stable. 
  \end{proof}
\begin{lemma}\cite{MR2463978}
 Suppose that $D(A) \subset D(C)$, $\rho(A) \neq \emptyset$, and that for some (and hence for all) $\mu \in \rho(A)$, the operator $(A- \mu)^{-1}B$ is bounded on $D(B)$ and $C(A-\mu)^{-1}$ is closed on $H$. Then $\mathcal A$ is closable iff $S_{2}(\mu)$ is closable for some (and hence for all) $\mu \in \rho(A)$. In this case, the closure $\overline{\mathcal {A}}$ is given by
\begin{equation}\label{equ 16}
    \overline{\mathcal A} = \mu + \begin{bmatrix}
    I & 0\\
    C(A-\mu)^{-1} & I
    \end{bmatrix} 
    \begin{bmatrix} 
    A-\mu & 0\\
    0 & \overline{S_{2}(\mu)}
    \end{bmatrix}
    \begin{bmatrix}
    I & \overline{(A- \mu)^{-1} B}\\
    0 & I
    \end{bmatrix}
\end{equation}
independently of $\mu \in \rho(A)$, that is,
$D(\overline{\mathcal A})= \biggl\{(x, y)^{T}\in H \oplus K: x + \overline{(A-\mu)^{-1}B}y \in D(A), y\in D(\overline{S_{2}(\mu)}) \biggr\}$, $\overline{\mathcal A} (x,y)^{T} = \bigl((A-\mu)(x + \overline{(A-\mu)^{-1}B}y)+ \mu x, C(x + \overline{(A-\mu)^{-1}B}y) + (\overline{S_{2}(\mu)} + \mu )y\bigr)^{T}$, where $S_{2}(\mu) = E-\mu- C(A-\mu)^{-1}B$ is called the  Schur complement of $\mathcal A$.
\end{lemma}
\begin{lemma}\cite{MR2463978}
 Suppose that $D(E) \subset D(B)$, $\rho(E) \neq \emptyset$, and that for some (and hence for all) $\mu \in \rho(E)$, the operator $(E- \mu)^{-1}C$ is bounded on $D(C)$ and $B(E - \mu)^{-1}$ is closed on $K$. Then $\mathcal A$ is closable iff $S_{1}(\mu)$ is closable for some (and hence for all) $\mu \in \rho(E)$. In this case, the closure $\overline{\mathcal {A}}$ is given by
\begin{equation}\label{equ 17}
    \overline{\mathcal A} = \mu + \begin{bmatrix}
    I & B(E-\mu)^{-1}\\
    0 & I
    \end{bmatrix} 
    \begin{bmatrix} 
    \overline{S_{1}(\mu)} & 0\\
    0 & E-\mu
    \end{bmatrix}
    \begin{bmatrix}
    I & 0\\
    \overline{(E- \mu)^{-1} C} & I
    \end{bmatrix}
\end{equation}
independently of $\mu \in \rho(E)$, that is,
$D(\overline{\mathcal A})= \biggl\{(x, y)^{T}\in H \oplus K:  \overline{(E-\mu)^{-1}C}x+y \in D(E), x\in D(\overline{S_{1}(\mu)}) \biggr\}$, $\overline{\mathcal A} (x,y)^{T} = \bigl((\overline{S_{1}(\mu)}+ \mu)x + B(\overline{(E-\mu)^{-1}C}x+y), (E-\mu)(\overline{(E-\mu)^{-1}C}x +y) +\mu y \bigr)^{T}$, where $S_{1}(\mu) = A-\mu- B(E-\mu)^{-1}C$ is called the Schur complement of $\mathcal A$.
\end{lemma}
\begin{theorem}\label{thm 9}
 Suppose that $\|Cx\| \leq a\|Ax\|$, for all $x \in D(A) \subset D(C), a \geq 0 $ and $0\in \rho(A)$. The operator $A^{-1}B$ is bounded in $D(B)$ with $S_{2}(0) = E - CA^{-1}B$ is closable. Then $R(\overline{S_{2}(0)})$ is closed iff $R(\overline{\mathcal A})$ is closed.
\end{theorem}
\begin{proof}
 Here, $\overline{\mathcal A} = \begin{bmatrix}

    I & 0\\
    CA^{-1} & I
    \end{bmatrix} 
    \begin{bmatrix} 
    A & 0\\
    0 & \overline{S_{2}(0)}
    \end{bmatrix}
    \begin{bmatrix}
    I & \overline{A^{-1} B}\\
    0 & I
    \end{bmatrix}$, where $S_2({0}) = E -CA^{-1}B$.
    and $\overline{\mathcal A} (u,v)^{T} = \bigl(Au + A\overline{A^{-1}B}v, Cu + C\overline{A^{-1}B}v + \overline{S_{2}(0)}v\bigr)^{T}$, for all $(u,v)^{T}\in D(\overline{\mathcal A})$.
    First, we will show that $\overline{R(\mathcal A)} \subset R(\overline{\mathcal A})$ when $R(\overline{S_{2}(0)})$ is closed. However, the reverse inclusion is obvious by the definition of the closable operator. Let  $(x,y)^{T} \in \overline{R(\mathcal A)}$. Then there exists a sequence $\{(x_{n}, y_{n})^{T}\}$ from $D(\mathcal A)= D(A) \oplus D(B)\cap D(E)$ such that $Ax_{n} + By_{n} \to x$ and $Cx_{n} + Ey_{n}\to y$ as $n\to \infty$. Hence $Cx_{n} + CA^{-1}By_{n} \to CA^{-1}x$, as $n\to \infty$. Thus $(E-CA^{-1}B)y_{n} = S_{2}(0) y_{n} \to y- CA^{-1}x$, as $n\to \infty$. Since $R(\overline{S_{2}(0)}) = \overline{R(S_{2}(0)})$, there exists  $w\in D(\overline{S_{2}(0)})$ such that $\overline{S_{2}(0)}(w) = y- CA^{-1}x$.
    Now, $\overline{\mathcal A} (A^{-1}x- \overline{A^{-1}B}w, w)^{T} = (x,y)^{T} \in R(\overline {\mathcal A})$. Therefore, $R(\overline{\mathcal A})$ is closed. 
    
  To prove the converse part, we will show that $\overline{R(S_{2}(0))} \subset R(\overline {S_{2}(0)})$, when $R(\overline{\mathcal A})$ is closed. Let  $y\in \overline{R(S_{2}(0))}$.  Then there exists a sequence $\{x_{n}\}$ from $D(S_{2}(0))$ such that $(E- CA^{-1}B)x_{n} \to y$, as $n\to \infty$.
Hence     $\overline{\mathcal A} (- \overline{A^{-1}B}x_{n}, x_{n})^{T} = (0, S_{2}(0)x_{n})^{T} \to (0,y)^{T}$, as $n\to \infty$. Since $R(\overline{\mathcal A})$ is closed, so there is an element $(u_{0}, v_{0})^{T} \in D(\overline{\mathcal A})$ such that $Au_{0}+ A\overline{A^{-1}B}v_{0} =0$ and $Cu_{0} + C\overline{A^{-1}B}v_{0} + \overline {S_{2}(0)}v_{0} =y$. From the inequality $\|Cx\| \leq a \|Ax\|$, we get that $\overline {S_{2}(0)}v_{0} =y \in R(\overline{S_{2}(0)})$. Therefore, $R(\overline{S_{2}(0)})$ is closed.
    
\end{proof}
\begin{corollary}
 Suppose that $\|Cx\| \leq a\|Ax\|$, for all $x \in D(A) \subset D(C), a \geq 0 $ and $0\in \rho(A)$. The operator $A^{-1}B$ is bounded in $D(B)$ with $S_{2}(0) = E - CA^{-1}B$ is a closable operator and $N(\overline{S_{2}(0)}) = \overline{N(S_{2}(0))}$. Then $S_{2}(0)$ is Hyers-Ulam stable iff $\overline{\mathcal A}$ is Hyers-Ulam stable.
\end{corollary}
\begin{theorem}
 Suppose that $\|Bx\| \leq b\|Ex\|$, for all $x \in D(E) \subset D(B), b \geq 0 $ and $0\in \rho(E)$. The operator $E^{-1}C$ is bounded in $D(C)$ with $S_{1}(0) = A - BE^{-1}C$ is closable. Then $R(\overline{S_{1}(0)})$ is closed iff $R(\overline{\mathcal A})$ is closed. Moreover, if  $N(\overline{S_{1}(0)})=\overline {N(S_{1}(0))}$, then $S_{1}(0)$ is Hyers-Ulam stable iff $\overline{\mathcal A}$ is Hyers-Ulam stable.
\end{theorem}
\begin{proof}
The proof is similar to  Theorem \ref{thm 9}.
\end{proof}
\begin{lemma}\cite{MR2463978}
 Suppose that $D(C) \subset D(A)$, $C$ is boundedly invertible, and that $C^{-1}E$ is bounded on $D(E)$ and $(A-\mu)C^{-1}$ is closed. Then $\mathcal A$ is closable (closed, respectively) iff $T_{2}(\mu)$ is closable (closed, respectively) for some (and hence for all) $\mu \in \mathbb {C}$. In this case, the closure $\overline{\mathcal A}$ is given by
\begin{equation}\label{equ 18}
    \overline{\mathcal A} = \mu + \begin{bmatrix}
        I & (A-\mu)C^{-1}\\
        0 & I
    \end{bmatrix}
    \begin{bmatrix}
    0 & \overline{T_{2}(\mu)}\\
    C & 0
    \end{bmatrix}
    \begin{bmatrix}
    I & \overline{C^{-1}(E-\mu)}\\
    0 & I
    \end{bmatrix}
\end{equation}
independently of $\mu\in \mathbb C$, that is,
$D(\overline {\mathcal A}) = \bigl\{(x,y)^{T}\in H \oplus K: x+ \overline{C^{-1}(E-\mu)}y \in D(C), y\in D(\overline{T_{2}(\mu)})\bigr \}$,
$\overline{\mathcal{A}} (x,y)^{T} = \bigl( (A-\mu)(x + \overline{C^{-1}(E- \mu)}y) + \mu x + \overline{T_{2}(\mu)}y, C(x + \overline{C^{-1}(E- \mu)}y) + \mu y \bigr )^{T}$, where $T_{2}(\mu) = B - (A-\mu)C^{-1}(E-\mu)$ is called the quadratic complement of $ {\mathcal A}$.

\end{lemma}
\begin{lemma}\cite{MR2463978}
 Suppose that $D(B) \subset D(E)$, $B$ is boundedly invertible, and that $B^{-1}A$ is bounded on $D(A)$ and $(E-\mu)B^{-1}$ is closed. Then $\mathcal A$ is closable (closed, respectively) iff $T_{1}(\mu)$ is closable (closed, respectively) for some (and hence for all) $\mu \in \mathbb {C}$.  In this case, the closure $\overline{\mathcal A}$ is given by
\begin{equation}\label{equ 19}
    \overline{\mathcal A} = \mu + \begin{bmatrix}
        I & 0\\
        (E- \mu)B^{-1} & I
    \end{bmatrix}
    \begin{bmatrix}
    0 & B\\
    \overline{T_{1}(\mu)} & 0
    \end{bmatrix}
    \begin{bmatrix}
    I & 0\\
    \overline{B^{-1}(A-\mu)} & I
    \end{bmatrix}
\end{equation}
independently of $\mu\in \mathbb C$, that is,
$D(\overline {\mathcal A}) = \bigl\{(x,y)^{T}\in H \oplus K: y+ \overline{B^{-1}(A-\mu)}x \in D(B), x\in D(\overline{T_{1}(\mu)})\bigr \}$,
$\overline{\mathcal{A}} (x,y)^{T} = \bigl( B(\overline{B^{-1}(A-\mu)}x +y) + \mu x, (E-\mu)(\overline{B^{-1}(A-\mu)}x +y) + \overline {T_{1}(\mu)}x +\mu y \bigr )^{T}$, where $T_{1}(\mu) = C - (E-\mu)B^{-1}(A-\mu)$ is called the quadratic complement of $ {\mathcal A}$.

\end{lemma}
\begin{theorem}\label{thm 10}
     Suppose that $\|Ax\| \leq a \|Cx\|$, for all $x\in D(C) \subset D(A)$ with $a \geq 0$ and $C$ is boundedly invertible.  The operator $C^{-1}E$ is bounded on $D(E)$ and $T_{2}(0)$ is closable. Then $R(\overline {T_{2}(0)})$ is closed iff $R(\overline{\mathcal A})$ is closed.
\end{theorem}
\begin{proof}
 $\mathcal A$ is closable and $\overline{\mathcal A}$ is given by
\begin{equation}\label{equ 20}
    \overline{\mathcal A} = \begin{bmatrix}
        I & AC^{-1}\\
        0 & I
    \end{bmatrix}
    \begin{bmatrix}
    0 & \overline{T_{2}(0)}\\
    C & 0
    \end{bmatrix}
    \begin{bmatrix}
    I & \overline{C^{-1}E}\\
    0 & I
    \end{bmatrix}
\end{equation}
where $D(\overline {\mathcal A}) = \bigl\{(u,v)^{T}\in H \oplus K: u+ \overline{C^{-1}E}v \in D(C), v\in D(\overline{T_{2}(0)})\bigr \}$,
$\overline{\mathcal{A}} (u,v)^{T} = \bigl( A(u + \overline{C^{-1}E}v) + \overline{T_{2}(0)}v, C(u + \overline{C^{-1}E}v)\bigr )^{T}$, and $T_{2}(0) = B - AC^{-1}E$.

First we claim that $\overline {R(\mathcal A)} \subset R(\overline{\mathcal A})$  when $R(\overline{T_{2}(0)})$ is closed. The reverse inclusion is true by the definition of the closablility of $\mathcal A$. Let  $(x,y)^{T} \in \overline {R(\mathcal A)}$. Then there exists a sequence $\{(x_{n}, y_{n})\}$ from $D(\mathcal A)= D(C) \oplus D(B)\cap D(E)$ such that $Ax_{n} + By_{n} \to x$ and $Cx_{n} + Ey_{n} \to y$, as $n\to \infty$. Hence  $Ax_{n} + AC^{-1}Ey_{n} \to AC^{-1}y$ as $n \to \infty$. Thus, $ T_{2}(0) y_{n} = (B-AC^{-1}E)y_{n} \to x-AC^{-1}y$, as $n\to \infty$. Since $R(\overline{T_{2}(0)})$ is closed, there is an element $w\in D(\overline {T_{2}(0)})$ such that $x- AC^{-1}y = \overline{T_{2}(0)}w$. As $\overline{\mathcal A} (C^{-1}y-\overline{C^{-1}E}w, w)^{T} = (x,y)^{T} \in R(\overline{\mathcal A})$, we get that $R(\overline{\mathcal A})$ is closed.

To prove the converse, we will show that  $\overline {R(T_{2}(0))} \subset R(\overline{T_{2}(0)})$ when $R(\overline{\mathcal A})$ is closed. Let $y\in  \overline {R(T_{2}(0))}$.  Then there exists a sequence $\{x_{n}\}$ from $D(T_{2}(0))$ such that $T_{2}(0)x_{n} \to y$, as $n\to \infty$. Hence $\overline{\mathcal A}(-C^{-1}Ex_{n}, x_{n})^{T} = (\overline{T_{2}(0)}x_{n}, 0)^{T} \to (y, 0)^{T}$, as $n \to \infty$. Since $R(\overline{\mathcal A})$ is closed, $(y,0)^{T} \in R(\overline{\mathcal A})$. There exists $(u_{0}, v_{0})^{T} \in D(\overline{\mathcal A})$ such that $A(u_{0} + \overline{C^{-1}E}v_{0}) + \overline{T_{2}(0)}v_{0} =y$ and $C(u_{0} + \overline{C^{-1}E}v_{0})= 0$. From the given inequality, we get that $\overline{T_{2}(0)}v_{0} = y \in R(\overline{T_{2}(0)})$. Therefore $R(\overline{T_{2}(0)})$ is closed.
\end{proof}
\begin{corollary}
 Suppose that $\|Ax\| \leq a \|Cx\|$, for all $x\in D(C) \subset D(A)$ with $a \geq 0$ and  $C$ is boundedly invertible.  The operator $C^{-1}E$ is bounded on $D(E)$ and  $T_{2}(0)$ is closable with $N(\overline {T_{2}(0)})= \overline{N(T_{2}(0))}$. Then ${T_{2}(0)}$ is Hyers-Ulam stable iff $\overline{\mathcal A}$ is Hyers-Ulam stable.
\end{corollary}
\begin{theorem}
 Suppose that $\|Ex\| \leq b \|Bx\|$, for all $x\in D(B) \subset D(E)$ with $b \geq 0$ and $B$ is boundedly invertible.  The operator $B^{-1}A$ is bounded on $D(A)$ and $T_{1}(0)$ is closable. Then $R(\overline {T_{1}(0)})$ is closed iff $R(\overline{\mathcal A})$ is closed.
\end{theorem}
\begin{proof}
    The proof is similar to the proof of Theorem \ref{thm 10}.
\end{proof}
\begin{corollary}
 Suppose that $\|Ex\| \leq b \|Bx\|$, for all $x\in D(B) \subset D(E)$ with $b \geq 0$ and $B$ is boundedly invertible.  The operator $B^{-1}A$ is bounded on $D(A)$ and  $T_{1}(0)$ is closable with $N(\overline {T_{1}(0)})= \overline{N(T_{1}(0))}$. Then ${T_{1}(0)}$ is Hyers-Ulam stable iff $\overline{\mathcal A}$ is Hyers-Ulam stable.
\end{corollary}

\section{An Illustrative Example of Closable (Non-closed) Hyers-Ulam Stable Operator}

\begin{proposition}\label{pro 1}
	Let $T$ denote the multiplication operator by $\phi \in \mathcal C(\mathbb R)$ on $L^{2}(\mathbb R)= L^{2}(\mathbb R, \Sigma, \mu)$ with domain $D(T)= \{f\in L^{2} (\mathbb R) : \phi f\in L^{2} (\mathbb R)\}$. Then $\mathcal{C}_{0}^{\infty}$ is a core for $T$. Here, measure $\mu$ is assumed to be a regular Borel measure. The space $\mathcal C(\mathbb R)$ denotes the set of all continuous functions on $\mathbb R$ and  $\mathcal{C}_{0}^{\infty}$  is the space of infinitely many times differentiable real-valued functions with compact support.
\end{proposition}
\begin{proof}
	Define, $L_{2,0} (\mathbb R)= \bigl\{ f\in L^{2}(\mathbb R) : \exists  K>0, |f(x)| \leq K$ a.e. in $\mathbb R$ and $f(x)=0 $ a.e. in $\{x\in {\mathbb R} : |x| >K\} \bigl \}$.  Let $f\in D(T)$.  Define $M_{n} = \{x\in \mathbb R: |f(x)| < n\}$. Then $M_{n}$ belongs to the Borel $\sigma$-algebra and $M_{n} \subset M_{n+1}$,  for all $n\in \mathbb N$.
	
	We consider a finite open interval $M$ such that $\mu(M) < \infty $,  $\mu(M \cap M_{n}) < \infty$ and $M \cap M_{n} \in \Sigma$. Now  define,
	\begin{equation}\label{equ 21}
		g_{n}(x)=
		\begin{cases}
			f(x)  &\text{if}  ~x\in M\cap M_{n}\\
			0 &\text{if} ~ x\in M\setminus (M\cap M_{n}).
		\end{cases}
	\end{equation}
	By using Lusin's theorem, for given $\frac{1}{n} >0$, there exists a compact set $K_{n}$ in $M\cap M_{n}$ such that  $h_{n} (x) = g_{n}(x) =f(x)$ on $K_{n}$ is continuous with $\mu((M_{n} \cap M )\setminus K_{n}) < \frac{1}{n}$
	and $h_{n}(x) =0$, for all $x\in M\setminus K_{n}$. We can consider $K_{n} \subset K_{n+1}$,  for all $n\in \mathbb N$. [If $K_{n} \not\subset K_{n+1}$ for some $n\in \mathbb N$, then we take union of two compact set $K_{n}$ and $K_{n+1}$ as $K_{n+1}$]. Then $h_{n}(x) \to f(x)$ as $n\to \infty$ a.e. in $M$ and $Th_{n}(x) \to Tf(x)$ as $n\to \infty$ a.e. in $M$. By dominated convergence theorem, $h_{n} \to f$ and $Th_{n} \to Tf$ as $n\to \infty$ in $L^{2}(M)$. Moreover, $h_{n} \in L_{2,0} (M)$, for all $n\in \mathbb N$. For given $\varepsilon >0$, there exists $n_{0}\in \mathbb N$ such that $\|h_{n_{0}} -f\|_{T} < \frac{\varepsilon}{4} $ where $h_{n_{0}} \in L_{2,0}(M)$.
	We know $\mathcal {C}_{0}^{\infty} (M)$ is dense in $L^{2}(M)$. 
	\begin{equation}\label{equ 22}
		L^{2}(M)= \overline {\mathcal {C}_{0}^{\infty} (M)} \supset L_{2,0} (M).
	\end{equation}
	From (\ref{equ 22}), there exists a sequence $\{P_{n}\}$ in $ {\mathcal C}_{0}^{\infty}(M)$ such that $P_{n} \to h_{n_{0}}$ as $n\to \infty$ in $L^{2}(M)$.  Hence for large  $n_{1}\in \mathbb N$, we have $\|P_{n_{1}} -h_{n_{0}}\| < \frac{\varepsilon}{8}$ and   $\|T P_{n_{1}} - Th_{n_{0}}\| < \frac{\varepsilon}{8}$ because $\phi \in \mathcal C(\mathbb R) $ is bounded in the open interval $M$, where $P_{n_{1}} \in \mathcal{C}_{0}^{\infty}(M)$. This implies $\|P_{n_{1}}- h_{n_{0}}\|_{T} < \frac{\varepsilon}{4}$. Hence, $\|P_{n_{1}} -f\|_{T} < \frac{\varepsilon}{2}$. It is obvious that $\mathcal{C}_{0}^{\infty}(M) \subset \mathcal{C}_{0}^{\infty}(\mathbb R)$. Thus $P_{n_{1}} \in \mathcal{C}_{0}^{\infty}(\mathbb R)$. Therefore, $\mathcal{C}_{0}^{\infty}(\mathbb R)$ is dense with respect to $T$-norm.  Thus $\mathcal{C}_{0}^{\infty}(\mathbb R)$ is core for $T$.
\end{proof}
\begin{example}
	We will show a non-closed operator to be Hyers-Ulam stable with the help of the Theorem \ref{thm 2}. Let $T_{\phi}$ be a multiplication operator on $L^{2}(M)$, where $L^{2}(M) = L^{2} (M, \Sigma, \mu)$, $\mu$ is a Borel regular measure, $M \subset \mathbb R$  and $\mu(M) < \infty$. Define $T_{\phi}$ as $T_{\phi}(f) =\phi f$, for $f\in D(T_{\phi}) = \{f\in L^{2}(M) : \phi f \in L^{2}(M), \phi \in \mathcal{C} (M) , |\phi(x)| \geq 1\}$. Then $T_{\phi}$ is closed and $T_{\phi} = T_{\phi}^{*} =T_{\overline {\phi}}$, where $\overline{\phi}$ is the complex conjugate  of $\phi$.  Moreover, $R(T_{\phi}) =R{(T_{\phi}}^{*}) = L^{2}(M)$ \cite{MR566954} and $N(T_{\phi}) = \overline {N(T_{\phi}|_{{\mathcal{C}}_{0}^{\infty}})} = \{0\}$, where ${T_{\phi}} |_{{\mathcal{C}}_{0}^{\infty}}$ is the restriction of $T_{\phi}$ in the domain ${\mathcal{C}}_{0}^{\infty}$. This shows that $T_{\phi}$ is Hyers-Ulam stable. By Theorem \ref{thm 2} and Proposition \ref{pro 1}, it can be shown that ${T_{\phi}} |_{{\mathcal{C}}_{0}^{\infty}} $ is Hyers-Ulam stable. But  ${T_{\phi}} |_{{\mathcal{C}}_{0}^{\infty}}$ is not closed. But it is closable because  ${T_{\phi}} |_{{\mathcal{C}}_{0}^{\infty}} \subsetneq \overline{{T_{\phi}} |_{{\mathcal{C}}_{0}^{\infty}}} = T_{\phi}$. 
\end{example}

 \begin{center}
	\textbf{Acknowledgements}
\end{center}

\noindent The present work of the second author was partially supported by National Board for Higher Mathematics (NBHM), Ministry of Atomic Energy, Government of India (Reference Number: 02011/12/2023/NBHM(R.P)/R\&D II/5947) and the work of the third author is supported by the grant from The Mohapatra Family Foundation.

\end{document}